\documentclass[preprint,10pt]{elsarticle}

\usepackage{natbib}        

\usepackage{amsfonts,amssymb,amstext,verbatim,amsmath,graphicx,color}
\usepackage{url}

\usepackage{amssymb}
\usepackage{latexsym}

\usepackage{graphics} %
\usepackage{epsfig} %
\usepackage{epstopdf} %

\newcommand{\ba}{\begin{array}}
\newcommand{\ea}{\end{array}}
\newcommand{\be}{\begin{equation}}
\newcommand{\ee}{\end{equation}}

\newtheorem{theorem}{\textbf{Theorem}}[section]
\newtheorem{lemma}[theorem]{\textbf{Lemma}}
\newtheorem{proposition}[theorem]{\textbf{Proposition}}

\newtheorem{assumption}[theorem]{\textbf{Assumption}}

\begin{document}
\begin{frontmatter}
\title{Performance Regulation of Event-Driven Dynamical Systems  Using
	Infinitesimal Perturbation Analysis}



\tnotetext[]{{Research  supported in
		part by  the NSF under Grant CNS-1239225.}  }

\author[First]{Y. Wardi}
\author[Second]{C. Seatzu}
\author[Third]{X. Chen}
\author[Fourth]{S. Yalamanchili}

\address[First]{School of Electrical and Computer Engineering,
	\\ Georgia Institute of Technology, Atlanta, Georgia, USA \\ (e-mail:
	ywardi@ece.gatech.edu).}
\address[Second]{Department of Electrical and Electronic
Engineering,\\ University of Cagliari, Italy \\ (e-mail:
seatzu @ diee.unica.it).}
\address[Third]{School of Electrical and Computer Engineering,
	\\ Georgia Institute of Technology, Atlanta, Georgia, USA \\ (e-mail:
	xchen318@gatech.edu).}
\address[Fourth]{School of Electrical and Computer Engineering,
	\\ Georgia Institute of Technology, Atlanta, Georgia, USA \\ (e-mail:
	sudha@ece.gatech.edu).}


\begin{abstract}
 This paper presents a performance-regulation method for a class of stochastic timed event-driven systems
aimed at output tracking of a given reference setpoint. The systems are either Discrete Event Dynamic Systems (DEDS)
such as queueing networks or Petri nets, or Hybrid Systems (HS) with time-driven dynamics and event-driven dynamics,
like fluid queues and hybrid Petri nets.
The regulator,
designed for simplicity and speed of computation, is comprised of  a single integrator having a variable gain to ensure
 effective tracking under time-varying plants.
 The gain's computation  is based on the Infinitesimal Perturbation Analysis (IPA) gradient of the plant function with
respect to the control variable, and the resultant tracking  can be quite robust with respect to modeling inaccuracies and gradient-estimation errors.
The proposed  technique is tested on examples taken from various application
areas and modeled with different formalisms, including queueing models, Petri-net model of a production-inventory
 control system, and a stochastic DEDS model of a multicore chip control.
Simulation results are presented in support of the proposed approach.
 \end{abstract}

\end{frontmatter}

\section{Introduction}

This paper describes a regulation technique for a class of dynamical systems,
designed for output tracking of a given setpoint reference.
The regulator consists of an integral control with a variable gain, computed on-line so as to enhance
the closed-loop system's stability margins and yield effective tracking. The gain-adjustment algorithm is based on
the derivative of the plant's output  with respect to  its input control, and therefore  the regulation technique
is suitable for systems  where such derivatives are readily computable in real time. This includes a class of
stochastic timed Discrete Event Dynamic Systems (DEDS) and Hybrid Systems (HD) where the derivative is computable by the
Infinitesimal Perturbation Analysis (IPA) sample-gradient technique.
 Our motivation is derived from the problem of
regulating instructions' throughput in multicore computer processors, and following an initial study of that problem in Ref. \cite{Almoosa12a}
we extend  the technique to a general class  of DEDS and HS.

The need for regulating instruction throughput at the hardware level in modern computer processors stems from
real-time applications where constant throughput
facilitates effective real-time task and thread (subprogram) scheduling, as well as  from multimedia
applications where  a fixed frame rate
must be maintained to avoid choppy video or audio. The design of effective regulators is challenging
 because of the lack of predictive analytical or prescriptive models,  and  unpredictable high-rate fluctuations of instructions-related
 switching activity factors at the cores.   For this reason, we believe, most of the published control techniques are ad hoc
 (see the survey in Ref. \cite{Lohn11}). A systematic control-theoretic  approach has been pursued in Refs.
 \cite{Brinkschulte09,Bauer10,Lohn11} which applied a PID controller and analyzed  the effects of
 proportional
 controls with fixed gains. Concerned with the  unpredictability
 and rapid changes in the thread-related activity    factors, Ref. \cite{Almoosa12a}  sought a controller with adaptive gain.
 Furthermore, it considered scenarios where measurements and computations in the control loop must be performed as
 quickly as possible, even at the expense of accuracy.  To this end it considered
 controlling the instruction throughput by a core's clock rate, and applied an integral controller
  whose real-time gain-adaptation algorithm is designed
 for stabilizing the closed-loop system and yielding   effective tracking convergence. The gain-adaptation algorithm is based on IPA
 as described in the sequel.

An abstract, discrete-time  configuration of the closed-loop system is shown in Figure 1,
where $n$ denotes the time-counter, $r$ is the setpoint reference, $u_n$ is the control input to the plant,   $y_n$ is the
resulting output, and $e_n:=r-y_n$ is the error signal. The system is single-input-single-output so that all
the quantities $u_n$, $y_n$, $e_n$ and $r$ are scalar.

Let $J:R\rightarrow R$ represent a performance function of the plant with respect to its input $u$, and assume that the function $J(u)$
is  differentiable. Given the $nth$ input variable $u_n$,  suppose that
the the plant's output $y_n$ provides an estimation of $J(u_n)$.
The controller that we consider has the form
\begin{equation}
u_n=u_{n-1}+A_ne_{n-1},
\end{equation}
and we recognize this as the discrete-time version of an integrator (summer) with a variable gain.
As mentioned earlier, the gain sequence $\{A_n\}$ is designed to enhance the stability margins of the closed-loop system and
reduce oscillations of the tracking algorithm while speeding up its convergence. As we shall see, one way to achieve that
is to have $A_{n}$ be defined as
\begin{equation}
A_n=\Big(J^{\prime}(u_{n-1})\Big)^{-1},
\end{equation}
with ``prime'' denoting derivative with respect to $u$.
However,  it may not be possible to compute the derivative term
$J^{\prime}(u_{n-1})$, and approximations have to be used. Denoting the approximation error by
$\phi_{n-1}$, the computed gain $A_n$ is defined as
\begin{equation}
A_n=\Big(J^{\prime}(u_{n-1})+\phi_{n-1}\Big)^{-1}.
\end{equation}

In the systems considered in this paper the  plant  represents average measurements taken from a physical system or a cyber system
over contiguous time-intervals
called {\it control cycles}. For example, suppose that the physical
 system is a continuous-time dynamical system with input $\upsilon(t)$ and output $\zeta(t)$, $t\geq 0$;
 its state variable is immaterial for the purpose of this discussion. Divide the  time axis into contiguous control cycles $C_n$, $n=1,2,\ldots$,  suppose that  the control input
is fixed during $C_{n}$ to a value $u_n:=v(t)$  $\forall t\in C_n$, and define $y_n$ by
\[
y_n\ :=\ \frac{1}{|C_n|}\int_{C_n}\zeta(t)dt,
\]
where $|C_n|$ is the duration of $C_n$. Alternatively, $y_n$ can represent average measurements taken from the output of a
discrete-time or discrete-event system. Generally we impose no restriction on the way the control
cycles are defined, they can be fixed a priori or determined by
counting events in a DEDS; we only require that
the input $u_n$ remains unchanged during $C_{n}$ and can be modified only when the next control cycle
begins.

\vspace{.2in}
\begin{figure}[h]
\centering
\includegraphics[width=0.85\textwidth]{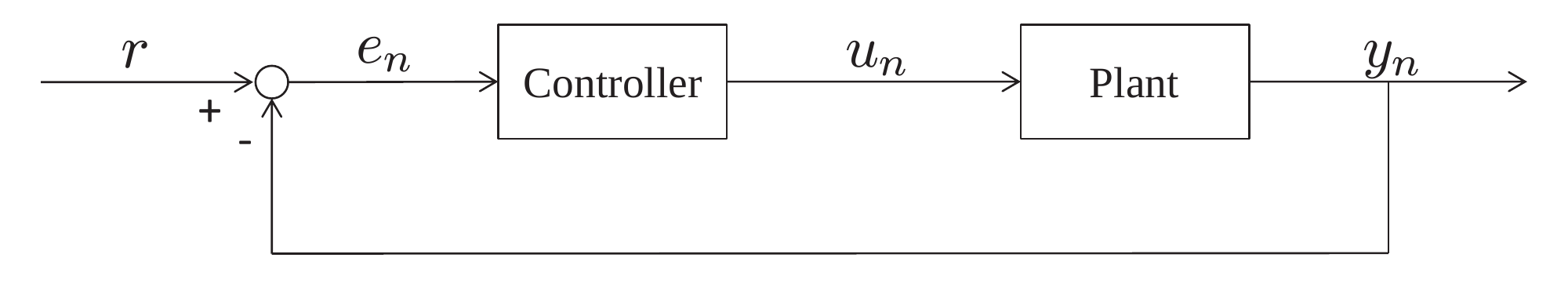}
{\small \caption{Basic regulation system}}
\end{figure}

Observe that Eq. (3) suggests that the computation of $A_{n}$ takes place during
the control cycle $C_{n-1}$. In fact, we assume that  the implementation of the control law takes place in the following temporal framework. Suppose that
the quantities  $u_{n-1}$, and $y_{n-1}$, $e_{n-1}$, and $A_{n}$ have been computed or measured by the starting time
of $C_{n}$.
Then $u_{n}$ is computed from Eq. (1) at the start of $C_{n}$ and we assume that this computation is immediate. During $C_{n}$,
 the plant produces $y_{n}$ from the applied input $u_{n}$ while $A_{n+1}$ is computed from
 Eq. (3), with the index $n+1$ instead of $n$. Finally, $e_{n}$ is computed at the end of $C_{n}$ from the equation
 \begin{equation}
 e_{n}=r-y_{n},
 \end{equation}
 and we assume that this computation is immediate.

 The plant's actions yielding  $y_{n}$ from $u_n$ during $C_n$ may represent a physical or cyber process or measurements thereof,
 and  the computation of $A_{n+1}$ is assumed to be carried out concurrently. Of a particular interest to us is the case where
 $J(u_n)$ is an expected-value performance function of a DEDS or HS,
 $y_{n}$ is an approximation thereof computed from a sample path of the system,  and the  term $J^{\prime}(u_n)+\phi_n$ in
 the Right-Hand Side (RHS) of Eq. (3) (with $n+1$)
 is computable by IPA.
 One of the main appealing features of IPA is the simplicity of its gradient (derivative)
 algorithms and efficiency of their computation.
 This, however, comes at the expense of accuracy. In particular,  in its principal  application area of queueing systems
 during its earlier development, IPA often yielded  statistically biased gradient estimators (see Refs. \cite{Ho91,Cassandras99}).
 To ameliorate this problem, Stochastic Flow Models (SFM) consisting of fluid queues (see Refs. \cite{Cassandras02,Sun04,Cassandras06,Panayiotou06}) and
 later extended to more general
 stochastic HS in Refs.   \cite{Cassandras10,Wardi10,Wardi13}, offer an alternative framework to queueing networks
 for the application  of IPA; in their setting   the IPA gradients typically are simpler and more accurate.
 Still approximations must be made either in the IPA algorithms or in the system's model when
 a stochastic HS is used as
 a
  modeling abstraction for a
DEDS. However, our overriding concern regarding the regulation's control law is that of simplicity and
 computational efficiency even if
 they come at the expense of accuracy.
 This is justified by  sensitivity-analysis  results, derived below,  showing that asymptotic tracking of the regulation
 scheme holds under substantial relative errors in the gradient estimation.

 The objective of this paper is to investigate the performance of our proposed tracking  technique on a number of
 DEDS and HS by using the IPA method for computing the integrator's adaptive gain in the loop. In this we leverage on the
 simplicity and low computational efforts required for the IPA derivatives. Furthermore,
 simulation experiments suggest that the regulation algorithm works well despite substantial errors in the
 gradient estimation, thus allowing  us to tilt the balance  between precision and low computing times towards
 fast computation at the expense of accuracy. It may be asked why we use an integral control and not a PI or PID controller, and it
 is pointed out that
 we have tested via simulation (not reported here) the
 addition of a proportional element to the integral control, and found no improvement. This is not surprising since, as indicated
 by the analysis in Section 2, the particular gain-adaptation of the integral controller stabilizes the system for a class of
 plant functions $J(u)$. In summary, the contributions of this paper are: 1). It  proposes the first general-purpose, systematic   performance-regulation technique for a class of timed DEDS and HS. 2). To-date, the main use of IPA has been in optimization, while this paper
 pursues a new kind of application, namely performance regulation. 3). It moves away from the traditional pursuit of unbiased IPA gradients, and instead searches for low-complexity approximations perhaps at the expense of accuracy or unbiasedness. We believe that these three
 points of novelty may open up a new  dimension in the research and applications   of IPA.

  The rest of the paper is organized as follows.
 Section 2  describes the regulation technique in an abstract setting. The following three sections present
 simulation results on three types of systems, and highlight the tradeoffs between simplicity and accuracy.
 Section 3 concerns a queue, Section 4 considers  a fluid Petri net production-control model, and
 Section 5 discusses a DEDS model of throughput  in computer processors. Finally,  Section 6 concludes the paper.\footnote{The second queueing example in Section 3, and the Petri-net example in Section 4 were presented in, and are part of  \cite{Seatzu14}.}

\section{Regulation Algorithm in an Abstract Setting}

The rationale behind the choice of $A_{n}$ in Eq. (2) (if its RHS  can be computed exactly) can be seen by considering the simple case where
the plant consists of a deterministic, time-invariant, memoryless nonlinearity, and hence the $u_n$~-~to~-~$y_n$ relation has the form
\begin{equation}
y_n=J(u_n).
\end{equation}
 Suppose that the function $J(u)$ is differentiable,
and denote its derivative by $J^{\prime}(u)$.
An application of the Newton-Raphson method for solving the
equation
\begin{equation}
r-J(u)=0
\end{equation}
results in the recursive equation
\begin{equation}
u_n=u_{n-1}+\frac{1}{J^{\prime}(u_{n-1})}e_{n-1}=u_{n-1}+A_ne_{n-1},
\end{equation}
where  $A_{n}:=\big(J^{\prime}(u_{n-1})\big)^{-1}$ is defined by (2). This yields
Eq. (1), and we discern that the control law, comprised of repeated recursive applications of Eqs.
$(2)\rightarrow(1)\rightarrow(5)\rightarrow(4)$, amounts to an implementation of the Newton-Raphson method.

An application of this control technique to a class of dynamic, time-varying systems is one of the main objectives of the present paper.
Its convergence is underscored by   established results on the sensitivity of the Newton--Raphson method with respect to variations in
function and gradient-evaluations \cite{Lancaster66}. Therefore, the  analysis  of the paper next will be presented
in an abstract setting of the Newton-Raphson method, and then related to the control configuration of Figure 1.

 Let $g:R\rightarrow R$ be a continuously-differentiable function, and consider
the problem of finding a root of the equation $g(u)=0$, $u\in R$. The basic step of the Newton-Raphson method is
\begin{equation}
u_{n}=u_{n-1}-\frac{1}{g^{\prime}(u_{n-1})}g(u_{n-1}),
\end{equation}
 and the algorithm is comprised of running this equation recursively
for $n=1,2,\ldots$ starting from an initial guess $u_0\in R$. Convergence of the algorithm can be characterized by
the equation
\begin{equation}
\lim_{n\rightarrow\infty}g(u_n)=0,
\end{equation}
which implies that every accumulation point of the sequence $\{u_n\}$, $\hat{u}$, satisfies the equation
$g(\hat{u})=0$.

Convergence of the Newton-Raphson method is well known under broad assumptions (see, e.g., Ref.
\cite{Ortega70}).
Ref. \cite{Lancaster66} investigated the case where the derivative term $g^{\prime}(u_{n-1})$ in
Eq. (8) is approximated rather than evaluated exactly, and showed that convergence is maintained under suitable
bounds on the errors. Taking it a step further, we consider the case where  errors arise in the function
evaluations $g(u_{n-1})$ as well. Then Eq. (9) no longer can be expected, but (see \cite{Almoosa12})
the limit
$\limsup_{n\rightarrow\infty}|g(u_{n})|$ is bounded from above by a term that depends on the magnitude of the errors
in a suitable sense.
Specifically, let
$\psi_{n-1}$ and $\phi_{n-1}$ denote additive  error terms in the computations of the function
$g(u_{n-1})$ and its derivative $g^{\prime}(u_{n-1})$, respectively,
so that Eq. (8) is transformed into
\begin{equation}
u_{n}=u_{n-1}-\frac{1}{g^{\prime}(u_{n-1})+\phi_{n-1}}\big(g(u_{n-1})+\psi_{n-1}\big).
\end{equation}
Define the relative errors
${\cal G}_{n-1}:=\frac{|\psi_{n-1}|}{|g(u_{-1})|}$ and ${\cal E}_{n-1}:=\frac{|\phi_{n-1}|}{|g^{\prime}(u_{n-1})|}$.
The following results, Lemma 2.2 and Proposition 2.3,
 are proved under the following assumption.\footnote{A weaker result than Lemma 2.2 and Proposition 2.3 was stated by Proposition 2 in
\cite{Almoosa12}, but its statement is incorrect.}
\begin{assumption}
The function $g(\cdot)$ is continuously differentiable.
\end{assumption}
Given a closed interval $I:=[u_1,u_2]\subset R$, define $|g^{\prime}|_{I,min}:=\min\{|g^{\prime}(u)|~:~u\in I\}$,
and
$|g^{\prime}|_{I,max}:=\max\{|g^{\prime}(u)|~:~u\in I\}$.
 Consider the algorithm comprised of recursive runs
of Eq. (10). For every $n=1,\ldots$, define
\[
m_{n}:=\min\{m>n~:~g(u_{m})g(u_{n})\geq 0\};
\]
in other words, if $g(u_{n})\geq 0$ then $m_{n}$ is the next integer $m>n$ such that $g(u_{m})\geq 0$,
and if $g(u_{n})\leq 0$ then $m_{n}$ is the next integer $m>n$ such that $g(u_{m})\leq 0$.

\begin{lemma}
For every   $M>1$  and $\beta\in(0,M^{-1})$  there exist $\alpha\in(0,1)$ and
$\theta\in(0,1)$ such that, for every closed interval $I$ where the function $g(u)$ has  the following three
properties:  (i) $g(\cdot)$ is either monotone nondecreasing throughout $I$ or monotone non-increasing throughout $I$; (ii) $g(\cdot)$ is either convex throughout $I$ or concave
throughout $I$; and (iii)
$\frac{|g^{\prime}|_{I,max}}{|g^{\prime}|_{I,min}}<M$, the following
 holds: if, for some $n=1,\ldots$, (a)  for every $j=n,\ldots,m_{n}$,
$u_{j}\in I$;    (b) for every $j=n,\ldots,m_{n}-1$,
${\cal E}_{j}<\alpha$; and   (c) for every $j=n,\ldots,m_{n}-1$, ${\cal G}_{j}<\beta$,
then
\begin{equation}
|g(u_{m_{n-1}})|<\theta|g(u_{n-1})|.
\end{equation}
\end{lemma}
\hfill$\Box$

The proof, tedious but  based on standard arguments from convex analysis,
is relegated to the appendix.

A few remarks are due.

\begin{enumerate}
\item
In situations where $g(u)$ is the expected-value performance function of a DEDS
 or HS, it may be impossible to  verify some of the assumptions underscoring Lemma 2.2, such as
the continuous differentiability  of $g(u)$, bounds on the relative errors ${\cal E}_{k-1}$ and
${\cal G}_{k-1}$, or  bounds on the terms $g^{\prime}_{I,max}$ and $g^{\prime}_{I,min}$. We point out that  analysis techniques
for their verifications have
been developed in the literature on IPA, mainly for convex sample performance functions  (see, e.g., \cite{Cassandras99, Cassandras06, Cassandras10}). However, in other situations these
assumptions  may have
to be stipulated or  justified by empirical evidence derived, for instance, from  simulation.
\item
Suppose that  the variable $u$ has to be constrained to a closed interval $I:=[u_{min},u_{max}]$ satisfying the conditions of Lemma 2.2.
To ensure that $u_{n}$, $n=1,\ldots$ are contained in $I$, it is possible to modify Eq. (10) by following it with a projection onto
$I$. Define the projection function $P_{I}:R\rightarrow I$ by
\[
P_{I}(u):=\left\{
\begin{array}{ll}
u, & {\rm if}\ u\in I\\
u_{min}, & {\rm if}\ u<u_{min}\\
u_{max}, & {\rm if}\ u>u_{max},
\end{array}
\right.
\]
and change Eq. (10)  to the following equation,
\begin{equation}
u_{n}=P_{I}\Big(u_{n-1}-\frac{1}{g^{\prime}(u_{n-1})+\phi_{n-1}}\big(g(u_{n-1})+\psi_{n-1}\big)\Big).
\end{equation}
If $I$ contains a point $\hat{u}$ such that $g(\hat{u})=0$, in addition to satisfying the conditions of the lemma, then it is readily seen
that using Eq. (12) instead of (10) does not weaken the statement of the lemma.
\end{enumerate}

\begin{proposition}
For every $\eta>0$,   $M>1$, and $\varepsilon>0$   there exist $\alpha\in(0,1)$ and $\delta>0$
 such that, for every closed, finite-length  interval $I$ such that (i) throughout $I$ the function
 $g(\cdot)$ is either monotone increasing or monotone decreasing, and either convex or concave,
  (ii) the set $\{u\in I~|~g(u)=0\}$ is nonempty, and
(iii) $|g^{\prime}|_{I,min}>\eta$ and $\frac{|g^{\prime}|_{I,max}}{|g^{\prime}|_{I,min}}<M$;
and for every sequence $\{u_n\}_{n=1}^{\infty}$
computed by a recursive application of Eq. (10) such that, for every $n=1,2,\ldots$,
(a) $u_n\in I$,
(b)
${\cal E}_{n}<\alpha$, and (c)  $|\psi_{n}|<\delta$,
 the following two inequalities are in force:
\begin{equation}
\limsup_{n\rightarrow\infty}|g(u_n)|<\varepsilon,
\end{equation}
and
\begin{equation}
\limsup_{n\rightarrow\infty}|g(u_n)+\psi_{n}|<\varepsilon.
\end{equation}
\end{proposition}
\hfill$\Box$

Remarks:
\begin{enumerate}
\item
Consider  the case where $g(u)$ is computed exactly and the errors are only in its derivative. Then every interval
$I$ and a sequence $\{u_{n}\}$ satisfying the conditions of the proposition also satisfy all the condition of Lemma 2.2, and hence,
Eq. (11) for all $n=1,\ldots$ implies (9).
\item
Notice that the conditions assumed in the proposition's statement include upper bounds on the {\it absolute} error of
the function's estimation, $|\psi_{n-1}|$, and on the {\it relative} error in the derivative, ${\cal E}_{n-1}$.
The reason for this discrepancy will be made clear in the proof of Proposition 2.3, where it will be shown  that
Eqs. (13) and (14) are due to such upper bounds.
\item
In light of Remark 2  following Lemma 2.2, if an interval $I$ satisfying the conditions of Proposition 2.3 is a constraint
set for the sequence $\{u_{n}\}$, then replacing Eq. 10) by (12) will not alter the assertions of the proposition expressed by Eqs. (13) and (14).
\end{enumerate}

Consider the case where $g(u)$ is an expected-value performance function on a stochastic dynamical system, and suppose that Eq. (10)
is run recursively in order to solve the equation $g(u)=0$. Suppose moreover that at the $nth$ iteration the function $g(u_{n-1})$
and its derivative $g^{\prime}(u_{n-1})$ are estimated by sample averages over
a control cycle $C_{n-1}$.  Generally, even under conditions of stochastic stability and arbitrarily long control cycles, it is not true
 that
$\{{\cal E}\}_{n-1}<\alpha$ and $|\psi_{n-1}|<\delta$ {\it for all} $n=1,\ldots$, as stipulated in the statement of Proposition 2.3. However practically,   by the inequality in Eq. (11),  with suitably-long
control cycles we expect the sequence $\{|g(u_n)|\}$ to decline towards  0 as in (13) except for sporadic jumps away from zero, and this
will be evident from the simulation results described in the next section. This is due to Lemma 2.2,
while Proposition 2.3 provides a unified result under more ideal conditions.

In the context of  the aforementioned control system, we note that the error in
 the derivative estimation is reflected by Eq. (3) rather than (2), and to account for the error in
 function evaluations, Eq. (5) is replaced by
 \begin{equation}
 y_n=J(u_n)+\psi_n.
 \end{equation}
Now the control loop is defined by a recursive application of Eqs. $(3)\rightarrow(1)\rightarrow(15)\rightarrow(4)$.
To translate the algorithm defined by Eq. (1) into this  control setting we define $g(u)=r-J(u)$. We  can also
apply it  to time-varying systems of the form $y_n=J_n(u_n)$, where convergence in the sense of Eqs. (13) -(14)
is expected
 if the systems vary slowly. Finally, we mention that from a practical standpoint the algorithm's convergence is ascertained
 from Eq. (14) rather than (13), since the observed quantity is $y_n=J(u_n)+\psi_n$ rather than $J(u_n)$.

The next section presents simulation examples of various DEDS and hybrid systems.
They include two noteworthy situations concerning queueing systems where IPA is biased.
The error term $\phi_{n-1}$ is due in one case to the use of   SFM as a basis for the IPA formula,
and in the other case, to a direct use of the DEDS-based IPA in spite of its bias. In both cases the regulation algorithm is shown to converge. While situations of the first case have stimulated an
interest in SFM as a means to circumvent the bias inherent in IPA,  the example of the second case
demonstrates that IPA can be successfully used even when it is biased.

 \section{Regulation of Average Delay and Loss Rate of a Single Queue}\label{sec:queues}
 This section illustrates the aforementioned regulation framework by applying it to two examples concerning,
 respectively,
 delay and loss rate in an M/D/1   queue.
 Both delay and loss  are controlled  by the service times. In the first example  the IPA derivative
 is unbiased and hence we use it to adjust the integrator's gain. In contrast, in the second example the IPA derivative is
 biased and hence  we apply a formula which is derived from an SFM (fluid-queue) approximation
 to the sample paths obtained from the discrete queue.  While this yields estimation errors due to modeling discrepancies, it
 guarantees the unbiasedness of IPA for the SFM and results in convergence of the regulation algorithm
 as applied to the discrete queue.

 \subsection{Average Delay}
 Consider an M/D/1/$\infty$ queue with a given arrival rate $\lambda$ and service times of $s$ time-units. Given positive integers $M$ and $N$, a control cycle consists
 of
 $M$ jobs,  and the regulation process is run for  $N$ cycles.  We set the queue to empty at the start of
 each control cycle.
 In the notation established  in Section 1,
 we define $u=s$;  for $n=1,\ldots,N$, $C_{n}$ is the $nth$ control cycle; and
 $y_{n}$ is the average delay  of jobs arriving during  $C_{n}$. Note that $y_n$ is not an expected-value
 function but rather a random function of $u$
  whose realization depends of the samples drawn during $C_{n}$. Let $J(u)$ denote the
 expected-value delay according to the stationary distribution. It is known that  $y_n$ is a strongly-consistent estimator
 of $J(u_n)$ in the following sense (see \cite{Ho91,Cassandras99}): for a given
 $n=1,2,\ldots$,
 \begin{equation}
 \lim_{M\rightarrow\infty}y_n=J(u_n),
 \end{equation}
 where $M$ is the length of the control cycle. Consequently the estimation error
 $\psi_n:=J(u_n)-y_n$ can be made smaller by choosing longer control cycles.

 The  IPA derivative, $\frac{\partial y_n}{\partial u_n}$,  is known
 to have the following form
 (see \cite{Ho91,Cassandras99}). For every $m=1,\ldots,M$, let $k_{m,n}$ denote the index (counter) of the job that started the busy
 period containing job $m$ during $C_{n}$. Then,
 \begin{equation}
 \frac{\partial y_n}{\partial u_n}=\frac{1}{M}\sum_{m=1}^{M}(m-k_{m,n}+1).
 \end{equation}
 In other words, the IPA derivative of the delay of job $m$ as a function of $u_n$ is the position of that job in its
 busy period, namely $m-k_{m,n}+1$. This IPA is unbiased and strongly consistent; see Refs. \cite{Ho91,Cassandras99}.
 Furthermore, by the Pollaczek-Khinchin formula,
  $J(u)$ is a continuously-differentiable function of $u$. It is obviously monotone increasing, and also  convex on $R^+$ since, by Eq. (17), $J^{\prime}(u)$ is monotone increasing as well. Therefore the conditions
 for Proposition 2.3 are satisfied on any closed interval contained in $R^+$,
 and we expect the simulation experiments to yield fast convergence of the output $y_n$
 to an $\varepsilon$-band around the target value $r$ (in the sense of Eq. (13)), where $\varepsilon$ can be made smaller by choosing longer control cycles.

 In the first simulation experiment we set the arrival rate to $\lambda=0.9$, and the target reference delay to $r=3.0$.
 The control cycles consist of $M=10,000$ jobs each, and we ran the control algorithm for $N=100$ cycles starting from
 the initial guess of $u_1=1.1$.  The results, shown in Figure 2, indicate an approach of $y_n$ to about a
 steady value in about 5 iterations.    Thereafter we notice fluctuations of $y_{n}$, and a closer view of the results, obtained in Figure 3 by focusing the graph on the range $n=6,\ldots,100$,  shows
 them to be in the range of $2.5 - 3.5$ except for a single exception at $n=69$.  These fluctuations do not seem to abate at larger $n$, and they likely are due
 to the variances of $y_{n}$ and its IPA derivative $\frac{\partial y_n}{\partial u_n}$. We point out that the average value of
 $y_{n}$, obtained over the range $n=6,\ldots,100$, is 3.031 (recall that the target value is 3.0).
 Correspondingly, the graph of the service times $u_n$, $n=1,\ldots,100$, is shown in Figure 4, and their mean over $n=6,\ldots,100$ is 0.913.

 \vspace{.2in}
 \begin{figure}[h]
 	\centering
 	\includegraphics[width=0.65\textwidth]{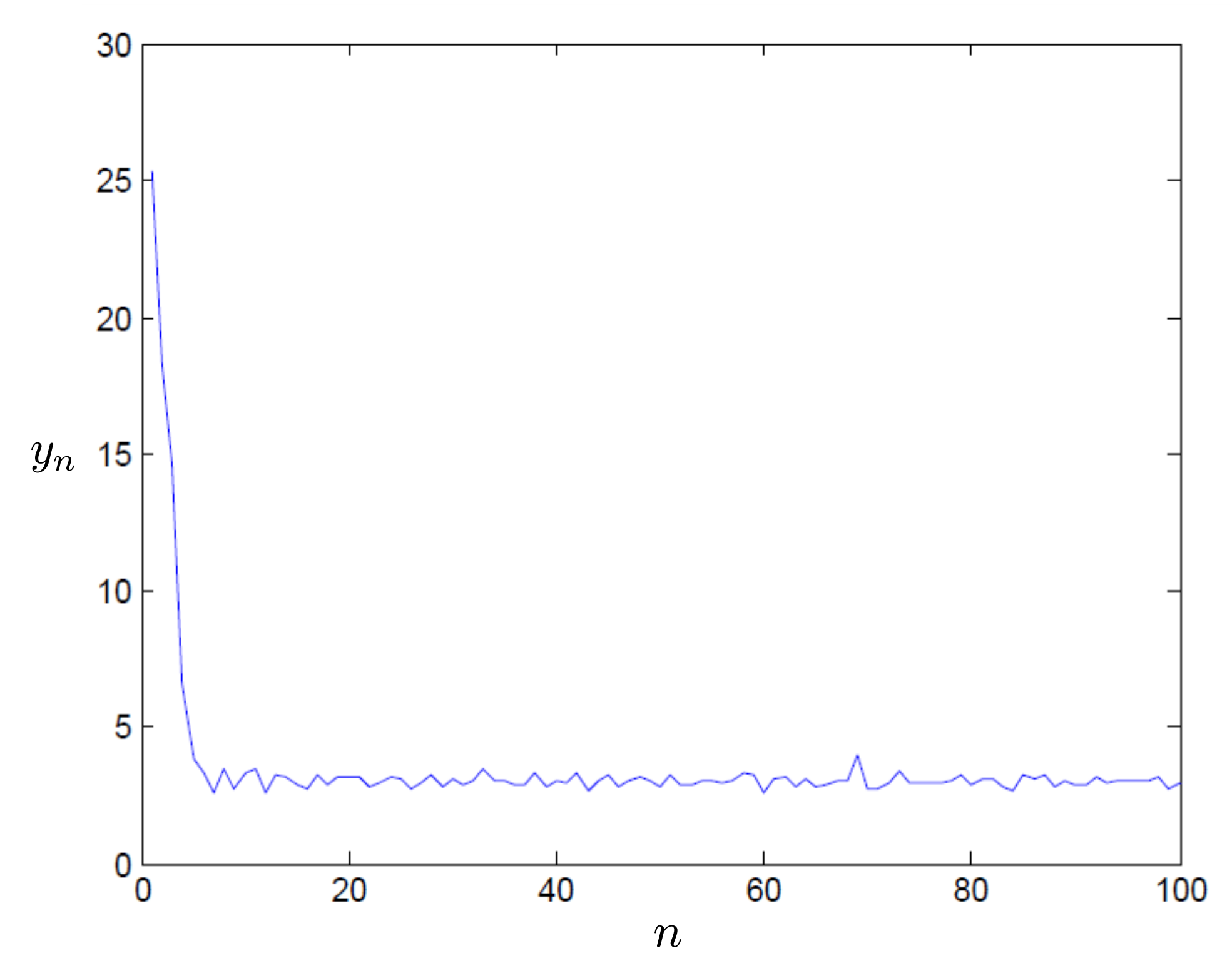}
 	{\small \caption{$M/D/1/\infty$ queue - delay.  $M=10,000$, $n=1,\ldots,100$.}}
 \end{figure}

 \vspace{.2in}
 \begin{figure}[h]
 	\centering
 	\includegraphics[width=0.65\textwidth]{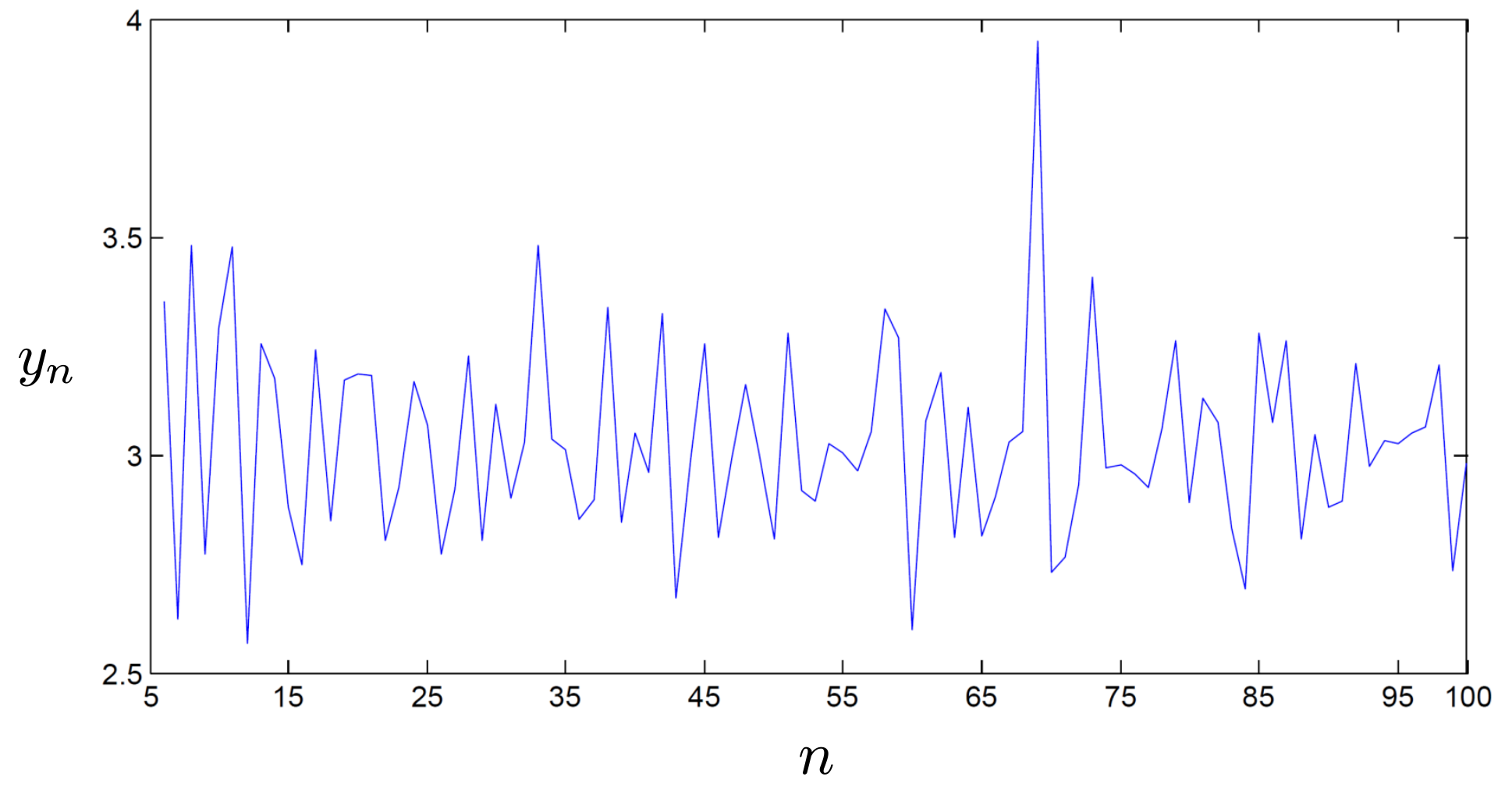}
 	{\small \caption{$M/D/1/\infty$ queue - delay, $n=6\ldots,\infty$.  $M=10,000$, $n=5,\ldots,100$.}}
 \end{figure}

 \vspace{.2in}
 \begin{figure}[h]
 	\centering
 	\includegraphics[width=0.65\textwidth]{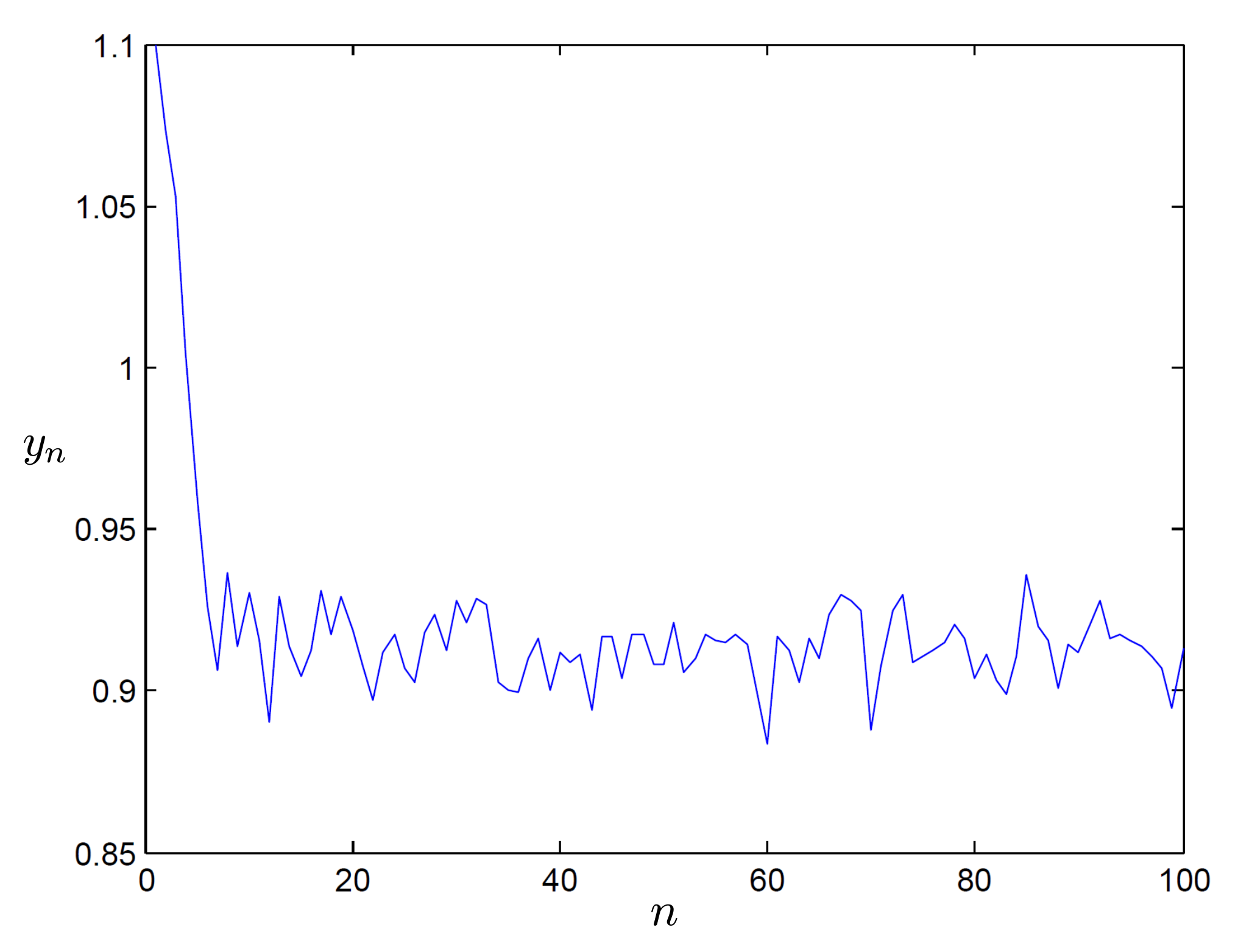}
 	{\small \caption{$M/D/1/\infty$ queue - service times.  $M=10,000$, $n=1,\ldots,100$.}}
 \end{figure}

 Lesser variance of $\{y_n\}$ can be  obtained from increasing the cycle length $M$ or reducing the reference target
 $y_{ref}$.
In the first case, where we set $M=100,000$ while keeping $y_{ref}$ at 3.00,   $y_n$,$, n=6,\ldots,100$,  fluctuated between
2.85 and 3.2 except for two values of $n$, and their mean was 3.008. In the second case, with $y_{{\rm ref}}=1.5$ and
$N=10,000$,
$y_{n}$ fluctuated mostly in the $[1.35-1.65]$-range with two exceptions, and its mean over $n=6,\ldots,100$ was
1.505.  In all three cases the regulation algorithm
 converged in about 5 iterations independently of the variances of the delays and their IPA derivatives; those seem to affect mostly the magnitude of the output fluctuations.

 In order to test the effects of changing the target reference during an experiment, we ran the regulation algorithm
 for 120 cycles starting with $u_1=1.1$; in the first 40 cycles $y_{ref}=3.0$, in the next 40 cycles $y_{ref}=4.5$, and in the last 40
 cycles, $y_{ref}=1.5$. The results are shown in Figure 5, and they indicate convergence to each value of $y_{ref}$ in about 5
 iterations. Furthermore, it is evident that the variations are larger for larger $y_{ref}$, and the reason is that
 the service times converge to larger values and hence the queue has larger traffic intensities. For
 $y_{ref}=3.0$ the mean of $y_{n}$ over $n=5,\ldots,40$ is 3.052,  for
 $y_{ref}=4.5$ the mean of $y_{n}$ over $n=45,\ldots,80$ is 4.653,  and for
 $y_{ref}=1.5$ the mean of $y_{n}$ over $n=85,\ldots,120$ is 1.504. The corresponding average gains $A_n$  obtained
 from the simulation were 0.063, 0.030, and 0.198. The respective  average service times were 0.918, 0.977, and 0.745. These values correspond to
 traffic intensities (product of $\lambda$ and the service time) of  0.826, 0.879, and 0.671, which explain the
 noticeable differences in variability. When we took $M=100,000$ the variability declined (as expected) and the obtained
 means are 3.041, 4.515, and 1.501, respectively.

 \vspace{.2in}
 \begin{figure}[h]
 	\centering
 	\includegraphics[width=0.65\textwidth]{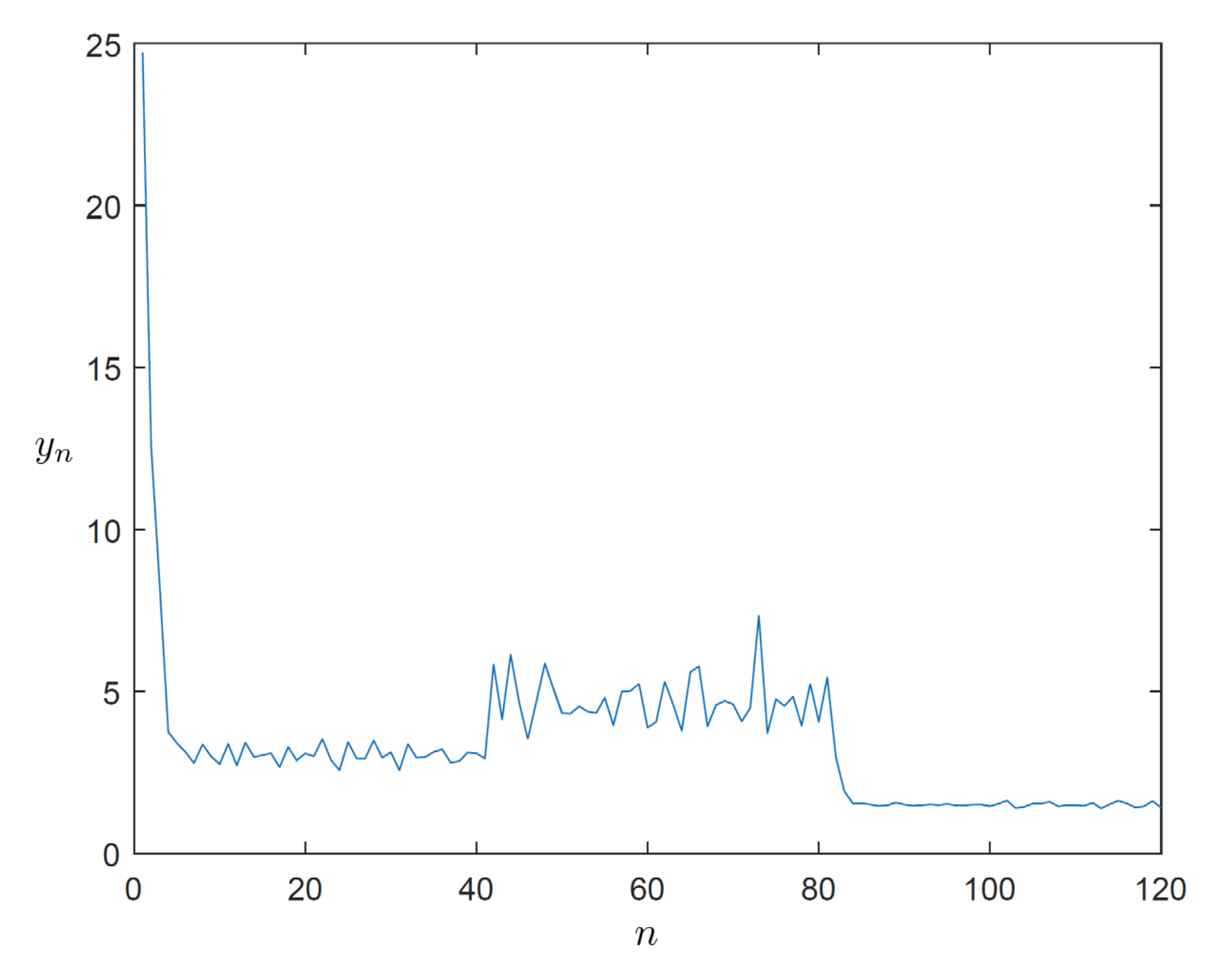}
 	{\small \caption{$M/D/1/\infty$ queue - delay, variable setpoint.  $M=10,000$, $n=1,\ldots,120$.}}
 \end{figure}

Lastly, we demonstrate the advantage of using the variable-gain controller over fixed-gain integral controls. To this
end we simulated the system described in the last paragraph with two fixed gains, and compared the results to those obtained
earlier from the variable-gain control. The constant gains are 0.030 and 0.198 which, as reported above, are the average
 gains driven by the variable-gain controller for the respective extreme setpoint references of 4.5 and 1.5.

 Figure 6 depicts the graph of $y_n$, $n=5,\ldots,120$, obtained from an application of the smaller
 fixed gain of 0.030 (the dashed curve), while the analogous graph obtained from the variable-gain
 controller  (the solid curve) is shown for the sake of comparison.
We discern slower response of the fixed-gain controller to changes in the setpoint,  and this is not
surprising since the variable-gain controller has larger gains  for the setpoints
 of 3.0 and 1.5.
 For  the larger fixed gain of 0.198, the graph of $y_n$ is shown in Figure 7.\footnote{This graph  could not be adequately shown in
 the same figure as the graphs in Figure 6  due to the large peaks of $y_n$.}   The controller performs well for the setpoint of 1.5, but
 it exhibits large oscillations suggesting instability for the setpoint of
 4.5, where the obtained gains of the variable-gain controller are smaller.
 These results are not surprising since, generally,
   tracking controllers with too-small gains may result in slow adjustment to variations in the target setpoint,
 while overly-large gains may result in oscillations and even instability. The variable-gain controller in this example
 seems to adjust well to changes in the setpoint and thus perform better than the considered fixed-gain controls.

 \vspace{.2in}
 \begin{figure}[h]
 	\centering
 	\includegraphics[width=0.65\textwidth]{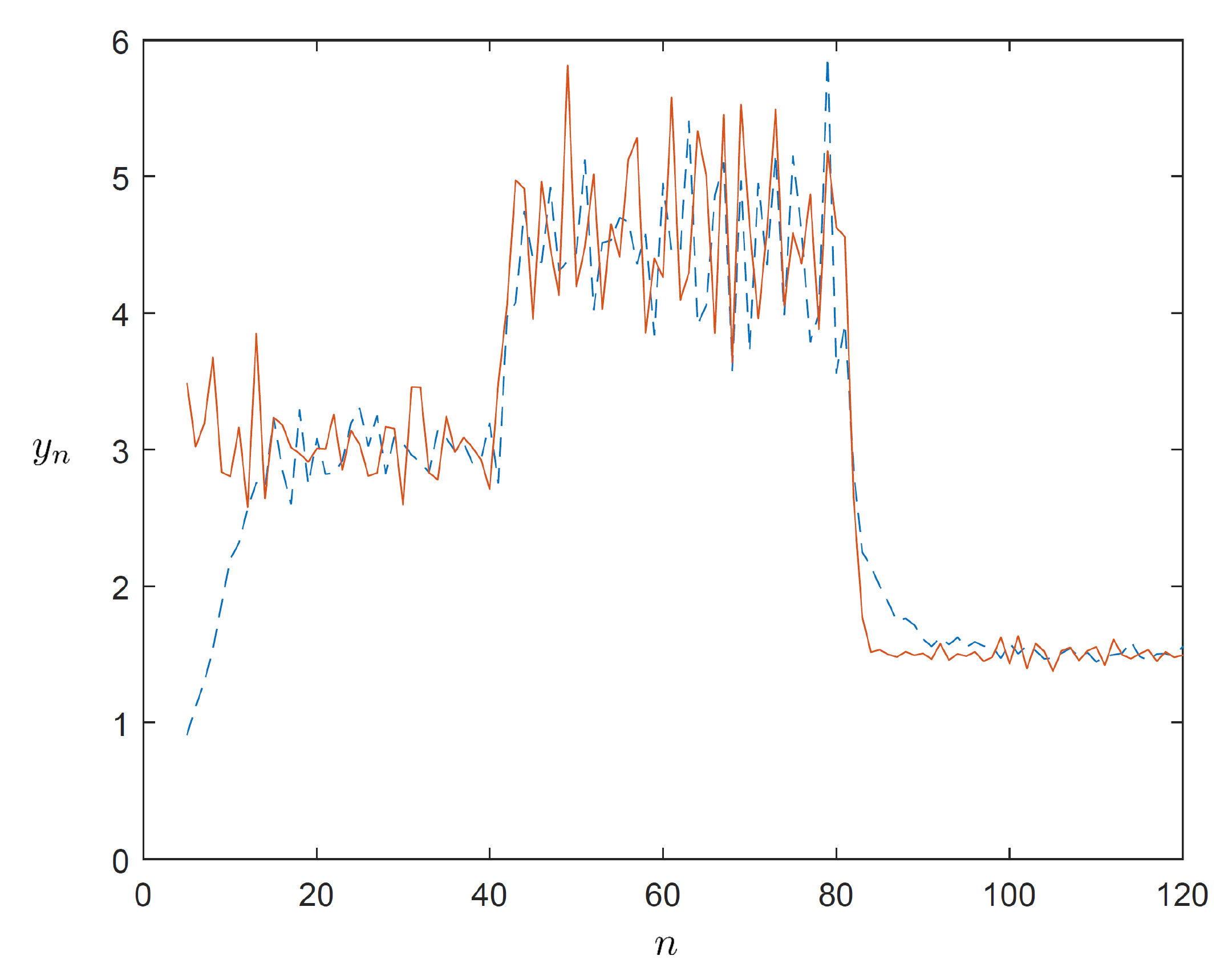}
 	{\small \caption{Comparison of a small, fixed-gain control to the variable-gain controller.}}
 \end{figure}

 \vspace{.2in}
 \begin{figure}[h]
 	\centering
 	\includegraphics[width=0.88\textwidth]{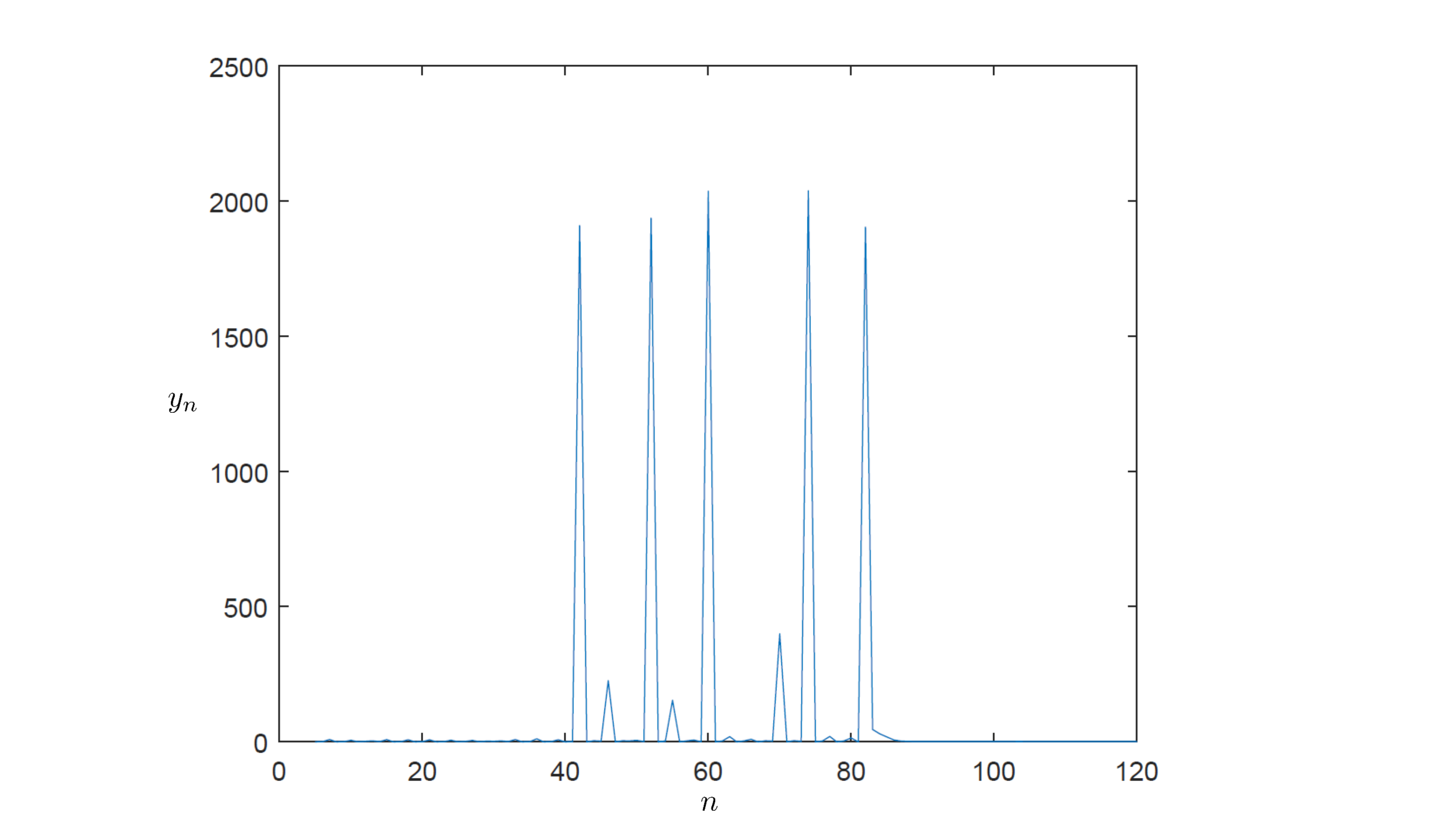}
 	{\small \caption{Oscillations due to a large, fixed-gain control.}}
 \end{figure}

 \subsection{Average loss-rate}

Consider  an M/D/1/k queue with a finite buffer, where jobs arriving at a full queue are
being discarded. Given an arrival rate, a buffer size, and a horizon period (cycle time), we aim at
controlling (regulating)    the mean loss rate   to a given reference by adjusting the service times.
 Accordingly, we denote the arrival rate by $\lambda$, the service times by $s$, the buffer size (including the holding
  place at the server) by $k$, and the horizon period by $t_{f}$.

  Let us divide the time-axis into consecutive control cycles, $C_{n}$, each of
  duration $t_{f}$ seconds. The control parameter is the service time, namely $u=s$, and during $C_{n}$ its
  value is denoted by
  $u_{n}$. The performance function of interest, $J(u)$, is the mean loss rate
  per control cycle. We approximate $J(u_n)$ by $y_n$, defined as  the number of  jobs
  discarded during $C_{n}$ divided by $t_{f}$.

  It is readily seen that the function $u_n\rightarrow y_n$ is piecewise
  constant and  almost surely no arrival would occur at the same time when another job enters or exits the server,
  and therefore, the IPA derivative  along a sample path  is $\frac{dy_n}{du_n}=0$. This does not provide an adequate
  approximation to $J^{\prime}(u_n)$, and hence an alternative approach to the estimation of
  $J^{\prime}(u_n)$ is needed.
  To this end we use a  fluid-queue SFM  as described
  in the next paragraph.

  Consider a fluid queue with a finite buffer, a time-varying inflow rate,  and a constant service
  rate. Suppose that the queue operates in a given time-interval $[0,t_{f}]$, where its instantaneous
  arrival rate, denoted by $\alpha(t)$, is a random process. Denote its service rate by $\beta$. Let $u:=\beta^{-1}$ be the control
  variable,
  and let $\gamma(u,t)$ denote the instantaneous overflow (spillover) rate due to the limited buffer.
  Let $L(u)$ denote the sample-based average loss rate per cycle as a function of the input service rate $u$, defined
   as
  \begin{equation}
  L(u)=\frac{1}{t_{f}}\int_{0}^{t_{f}}\gamma(u,t)dt.
  \end{equation}
  Observe that $L(u)$ provides an approximation to $J(u)$ to the extent that the fluid
  queue serves to approximate the discrete queue. The main difference between $J(u)$ and $L(u)$ is that $J(u)$ is the
  mean loss rate of the discrete queue, while $L(u)$ is the sample time-average loss rate per cycle of the continuous queue.

  Reference \cite{Cassandras02} showed that the IPA derivative $L^{\prime}(u)$ is unbiased, and
  it is computable by a simple, model-free  formula (listed below) that can act on the sample paths of the discrete queue.
  Furthermore, it can serve as an approximation to $J^{\prime}(u)$.
  Therefore we  implement the regulation algorithm in the following way.
  During cycle $C_{n-1}$, the service time is $u_{n-1}$,   $y_{n-1}$  is the resulting
  sample-average loss rate of the discrete queue, and $\psi_{n-1}:=y_{n-1}-J(u_{n-1})$. The IPA derivative  $L^{\prime}(u_{n-1})$,
  specified below,
   is the sample derivative of the average loss rate of the continuous queue as defined by Eq. (18),
   and $\phi_{n-1}:=L^{\prime}(u_{n-1})-J^{\prime}(u_{n-1})$.
  Eq. (3) becomes
  \begin{equation}
  A_{n}:=\big(L^{\prime}(u_{n-1})\big)^{-1}.
  \end{equation}
 In this  we apply
  the IPA derivative-formula,    obtained from an analysis of the SFM,  to the sample
  path of the discrete queue during $C_{n-1}$. In contrast, the plant's output $y_{n-1}$ that is used in the control loop via
  Eq. (1) with $e_{n-1}:=r-y_{n-1}$ corresponds to the discrete queue  since it represents the
  ``real'' system.
  The effectiveness  of the resulting regulation algorithm is related to the  quality of the approximation
  of $J(u_n)$ and $J^{\prime}(u_n)$ by $y_n$ and $L^{\prime}(u_n)$,
  respectively. Now Ref. \cite{Cassandras02} showed that the function $L(u)$ is Lipschitz continuous in
  $u$, w.p.1. Furthermore, by Eq. (20), below, $L^{\prime}(u)\geq 0$, and $L(u)$ and $L^{\prime}(u)$ are monotone increasing and hence
  the function $L(u)$ is convex w.p.1. Therefore we believe that $J(u)$ is convex as well although we do not know of a proof,
  and in this case the assumptions of  Lemma 2.2 and Proposition 2.3 are in force. In any event, the effectiveness of the regulation
  algorithm will be tested by simulation.

  The  IPA derivative $L^{\prime}(u)$ has the following form (see \cite{Cassandras02}).
   Suppose that there are $Q$ lossy busy periods during the horizon interval $[0,t_{f}]$, indexed
  by $q=1,\ldots,Q$ in increasing order (a busy period is {\it lossy} if  any positive amount of overflow
  is incurred throughout its duration).
  For the $qth$ lossy busy period, let $u_{q}$ be the first time loss occurs during it, and let $v_{q}$ be its end point; in other words,
  $u_{q}$ is the first time in that busy period when the buffer becomes full, and $v_{q}$ is the next time the buffer becomes empty. Then, under mild assumptions \cite{Cassandras02},
  \begin{equation}
  L^{\prime}(u)=\frac{1}{t_{f}}u^2\sum_{q=1}^{Q}\big(v_{q}-u_{q}\big).
  \end{equation}

  With $A_{n}$ defined by (19), we ran a simulation with the following parameters: $y_{{\rm ref}}=0.1$, $t_{f}=4,000$,
  $k=3$, $\lambda=0.9$, the initial parameter-value was $u_{1}=1.5$, and the number of cycles was $N=100$. The resulting
  graph of $y_{n}$, $n=1,\ldots,100$ is shown in Figure 8, where we notice convergence of the tracking algorithm in 3 iterations,
  to a band around the target value of 0.1. Within this band $y_{n}$ fluctuates between 0.08 and 0.125, except for a single value of $n$ where $y_{n}=0.139$. The mean of $y_{n}$ in the range $n=5,\ldots,100$ was 0.1005. To reduce the variability we ran the simulation for
  $t_{f}=20,000$, and the results, shown in Figure 9, exhibit an equally-fast convergence of the regulation algorithm
  with fluctuations in the range of $[0.091,0.114]$, and  with mean (over $n=5,\ldots,100$) of 0.1000.

\vspace{.2in}
\begin{figure}[h]
\centering
\includegraphics[width=0.65\textwidth]{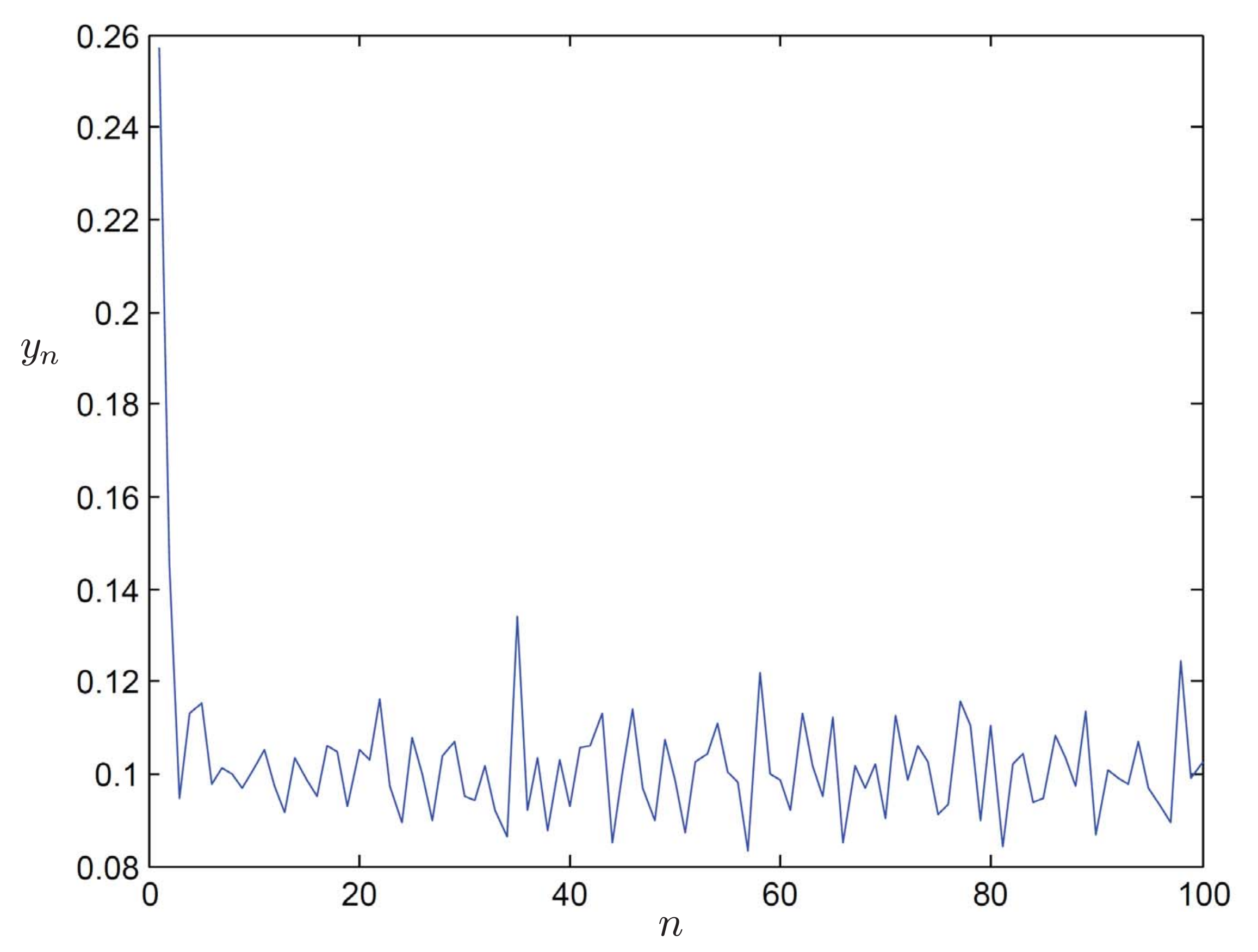}
{\small \caption{M/D/1/k queue - loss.  $t_{f}=4,000$, $n=1,\ldots,100$.}}
\end{figure}

\vspace{.2in}
\begin{figure}[h]
\centering
\includegraphics[width=0.65\textwidth]{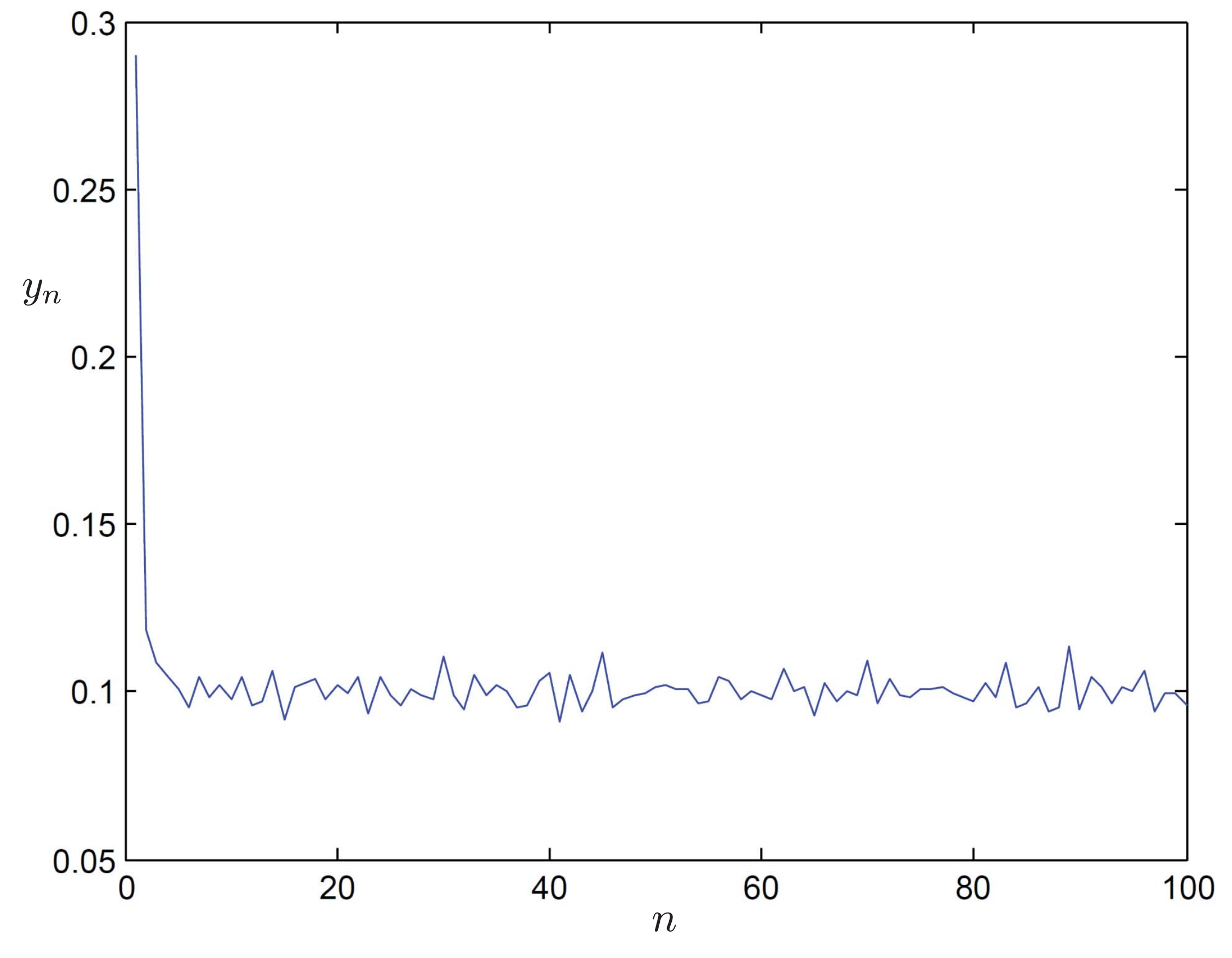}
{\small \caption{M/D/1/k queue - loss.  $t_{f}=20,000$, $n=1,\ldots,100$.}}
\end{figure}

\section{Application to a production system modeled by a Petri net}\label{sec:PN}

The IPA technique recently has been extended from fluid queueing networks
to a class of continuous Petri nets \cite{Xie02,Wardi13b}. References \cite{Wardi13b,Seatzu13} applied the results to an optimization example
of balancing part-inventories with product backorders in a single-stage manufacturing system,  and tested the
application of IPA in conjunction with a stochastic approximation algorithm. This section uses the same example to test
our approach to regulation.

The considered manufacturing system  consists of a machine that produces a sequence of single-type
 products.
The production schedule is driven by products' orders while  parts' inventories are maintained as safety stocks.
To make a product, the
system must have an available  part and a product order; parts without orders accumulate in the form of
inventories, while orders without parts result in cumulative backorders. Naturally
both excessive inventories and backorders are undesirable, and References
\cite{Wardi13b,Seatzu13} devise an IPA-based algorithm for optimally balancing
them. The underlying model for the algorithm is comprised of the
continuous (fluid) Petri net (event graph) shown in Figure 10.

\vspace{.2in}
\begin{figure}[h]
\centering
\includegraphics[width=0.65\textwidth]{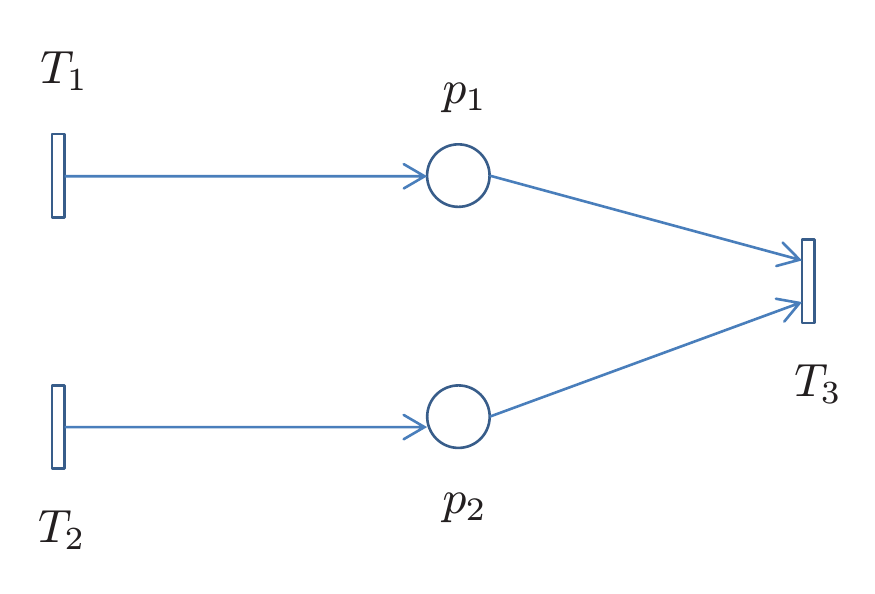}
{\small \caption{Basic  event graph}}
\end{figure}

Continuous Petri nets are hybrid Petri nets where the flow of fluid tokens through transitions is represented by
piecewise-continuous rate processes; see, e.g., \cite{Silva04} for comprehensive presentations.
With regard to our system shown in Figure 10,
 transitions
$T_{1}$, $T_{2}$, and $T_{3}$ represent, respectively, the processes of product-orders,
 parts' arrivals,
and the machine's operation. Each transition $T_{i}$, $i=1,2,3$, is
characterized by a maximum fluid-flow rate $V_{i}(t)>0$, which acts as an upper bound on its actual flow rate,
denoted by $v_{i}(t)$. The places $p_{1}$ and $p_{2}$ are used for holding fluid, and at time $t$ the amount
of stored fluid is denoted by $m_{1}(t)$ and $m_{2}(t)$, respectively. The processes $\{V_{i}(t)\}$, $\{v_{i}(t)\}$
($i=1,2,3$) and $\{m_{j}(t)\}$ ($j=1,2$) can be viewed as random processes defined over a common
probability space $(\Omega,{\cal F},P)$.

The dynamics of the system are described by the following three equations relating the above processes: For $i=1,2$,
\begin{equation}
v_{i}(t)=V_{i}(t).
\end{equation}
 For $i=3$, define $\varepsilon_{3}(t):=\big\{j\in\{1,2\}:m_{j}(t)=0\big\}$; then
\begin{equation}
v_{3}(t)=\left\{
\begin{array}{ll}
V_{3}(t), & {\rm if}\ \varepsilon_{3}(t)=\emptyset\\
min\big(v_{i}(t):i\in\varepsilon_{3}(t)\big), & {\rm if}\ \varepsilon_{3}(t)\neq\emptyset.
\end{array}
\right.
\end{equation}
As for $m_{j}(t)$, $j=1,2$, we have that
\begin{equation}
\dot{m}_{j}(t)=v_{j}(t)-v_{3}(t).
\end{equation}
In the forthcoming discussion we assume  that the system evolves in a given time-interval $[0,t_{f}]$ with given
initial conditions $m_{1}(0)$ and $m_{2}(0)$.

In the typical case where the three processes $\{V_{i}(t)\}$ are exogenous, the other network processes,
$\{v_{i}(t)\}$ and $\{m_{j}(t)\}$,
 are defined
in their terms via Eqs. (21)-(23). In other situations some of the processes $\{V_{i}(t)\}$ are exogenous
while others are controlled, and in this case the equations describing the controls together with (21)-(23)
define all of the network processes.
In the example considered in \cite{Wardi13b,Seatzu13} the processes $\{V_{1}(t)\}$ and $\{V_{3}(t)\}$ are
exogenous while $\{V_{2}(t)\}$ is controlled. Specifically, product orders are assumed to arrive in batches, and hence
\begin{equation}
V_{1}(t)=\sum_{k\geq 1}\alpha_{k}\delta(t-s_{k}),
\end{equation}
where $\delta(\cdot)$ is the Dirac delta function, $s_{k}$, $k=1,2,\ldots$,  are the arrival times, and $\alpha_{k}$ represent the quantities of the orders. The machine is assumed
to have deterministic service times, and hence $V_{3}(t)=V_{3}$ for a given $V_{3}>0$. The parts' arrival rates are controlled by the backorders via a threshold in the following fashion: $V_{2}(t)$ has a given low value  if the backorder levels  are below the threshold,
and a given higher value if the backorder levels are above the threshold. Formally, given a threshold $\rho>0$,
and given constants $V_{2,1}\geq 0$ and
$V_{2,2}>V_{2,1}$, $V_{2}(t)$ is defined via
\begin{equation}
V_{2}(t)=\left\{
\begin{array}{ll}
V_{2,1}, & {\rm if}\ m_{1}(t)\leq\rho\\
V_{2,2}, & {\rm if}\ m_{1}(t)> \rho.
\end{array}
\right.
\end{equation}
We assumed that $V_{2,1}\leq V_{3}\leq V_{2,2}$. It is obvious that Equations (21)-(23) and (25) have a unique joint solution for
every set of initial conditions.

Now let us consider the threshold $\rho$ as the control parameter and hence denote it by $u=\rho$.
 Then the processes
$\{V_{2}(t)\}$, $\{v_{i}(t)\}$, $i=2,3$, and $\{m_{j}(t)\}$, $j=1,2$ are functions of $u$ as well,
and hence are denoted by $\{V_{2}(u,t)\}$, $\{v_{i}(u,t)\}$, and $\{m_{j}(u,t)\}$, respectively. Assume that
a particular value of $u$ remains fixed throughout the evolution of the system in a given interval
$[0,t_{f}]$. Consider the sample performance
function $L(u)$ defined as
\begin{equation}
L(u)=\frac{1}{t_{f}}\int_{0}^{t_{f}}m_{2}(u,t)dt,
\end{equation}
for a given distribution of the initial conditions $m_{j}(0)$, $j=1,2$, and let $J(u):=E\big(L(u)\big)$
denote its expected value.
References \cite{Wardi13b,Seatzu13} considered the problem of minimizing
 a weighted sum of
$J(u)$ and the expected-value of the average workload at $p_{1}$. Here we use the same system
except that we perform regulation of $L(u)$ rather than optimization.

To put it all in the setting described
in Section 1, we define a control cycle to consist of $t_f$ time units, $u_n$ denotes the input during the $nth$ control
cycle $C_n$,
and $y_n:=L(u_n)$ as defined by Eq. (26). Therefore, for every $n=1,\ldots$,
$y_n\rightarrow J(u_n)$ as $t_f\rightarrow\infty$ w.p.1. Moreover,
Refs.  \cite{Wardi13b,Seatzu13} prove that the IPA derivative $L^{\prime}(u)$ is unbiased, hence
$\frac{\partial y_n}{\partial u_n}\rightarrow_{t_f\rightarrow\infty}J^{\prime}(u_n)$ w.p.1. Consequently, the error terms
$\psi_n$ and $\phi_n$ can be  reduced by taking  longer cycle times $t_f$. Regarding monotonicity and convexity of $J(u)$,
larger threshold $u$ results in smaller inventories and hence $J(u)$ is monotone non-increasing. However, we have
no way of proving convexity or concavity of $J(u)$, hence we put the regulation algorithm to the test
by simulation.

In the considered example, $V_{3}=6$, $V_{2,1}=2.15$, and $V_{2,2}=6$; these numbers are taken from
\cite{Wardi13b,Seatzu13}. The product-orders process $\{V_{1}(t)\}$, defined
by
(24), consists of equally-spaced arrivals every 50 seconds (deterministic), and each arrival
brings in an amount of fluid that is uniformly distributed in the $[30,70]$-range. The reference value to be tracked is
$J_{{\rm ref}}=758.70$, and it is the computed value of $J$ obtained for the aforementioned optimization problem in
\cite{Seatzu13}.
The IPA derivative $L^{\prime}(\theta)$ is computable via a recursive algorithm
constructed according to  the event-calculus framework defined in \cite{Cassandras10,Wardi13b};  a detailed presentation thereof
can be found in \cite{Seatzu13}.

We ran the regulation algorithm for 100 control cycles with $t_f=1,000$. The
results of two typical runs are shown in Figure 11 for two respective values of the initial control parameter,
$u_{1}=35$ and $u_{1}=15$. The corresponding  graphs of $y_{n}:=L_{n}(u_{n})$ are plotted
by the dashed curve and solid curve, and both indicate convergence to a band around $y_{{\rm ref}}=758.70$ after 3 iterations.
This band has
a maximum range of 47.17, and the averages of the outputs $y_{n}$, taken over $n=20,\ldots,100$, are 758.55 for the dashed plot,
and 758.73 for the solid plot. Although these results indicate convergence of the outputs' average to $y_{{\rm ref}}$, variations
in the output values are discernable.  These variations, as well as those in the IPA derivatives, yield
 fluctuations in the values of $u_{n}$, $n=1,2,\ldots$, as can be seen in Figure 12. We believe that the
 major cause of these variations is in the variance of $L(u)$ and   not inherent in the regulation
 algorithm. To test this point, we ran 100 independent simulations of the system at the fixed value of
 $u=24.8$, which is close to the
 average of $u_{n}$, $n=20,\ldots,100$ obtained by the regulation algorithm. The results,
 $\big\{L(u)\big\}_{i}$, shown in Figure 13,
 indicate a persistent variation with a maximum range of 42, which is comparable to the range obtained from Figure 11.

 To further test the convergence rate of the regulation algorithm we chose more extreme starting values of the control
 parameter, namely $u_{1}=45$ and $u_{1}=5$. The corresponding values of $y_{1}$ are 635 and 883, which are more extreme
 than those obtained in Figure 11, but nonetheless the regulation algorithm converged to a similar band around $y_{{\rm ref}}$
 in 3 iterations.

\vspace{.1in}
\begin{figure}[h]
\centering
\includegraphics[width=0.65\textwidth]{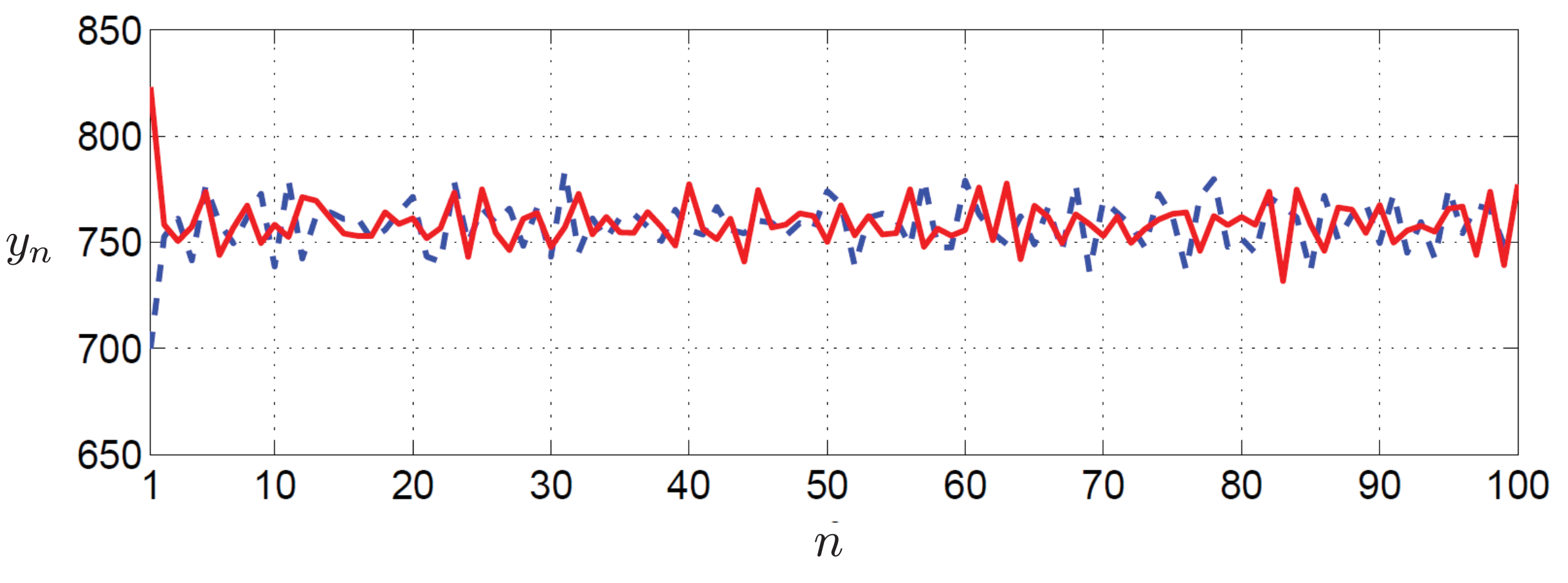}
{\small \caption{Petri net - inventory. $t_{f}=1,000$, $n=1,\ldots,100$}}
\end{figure}

\vspace{.1in}
\begin{figure}[h]
\centering
\includegraphics[width=0.65\textwidth]{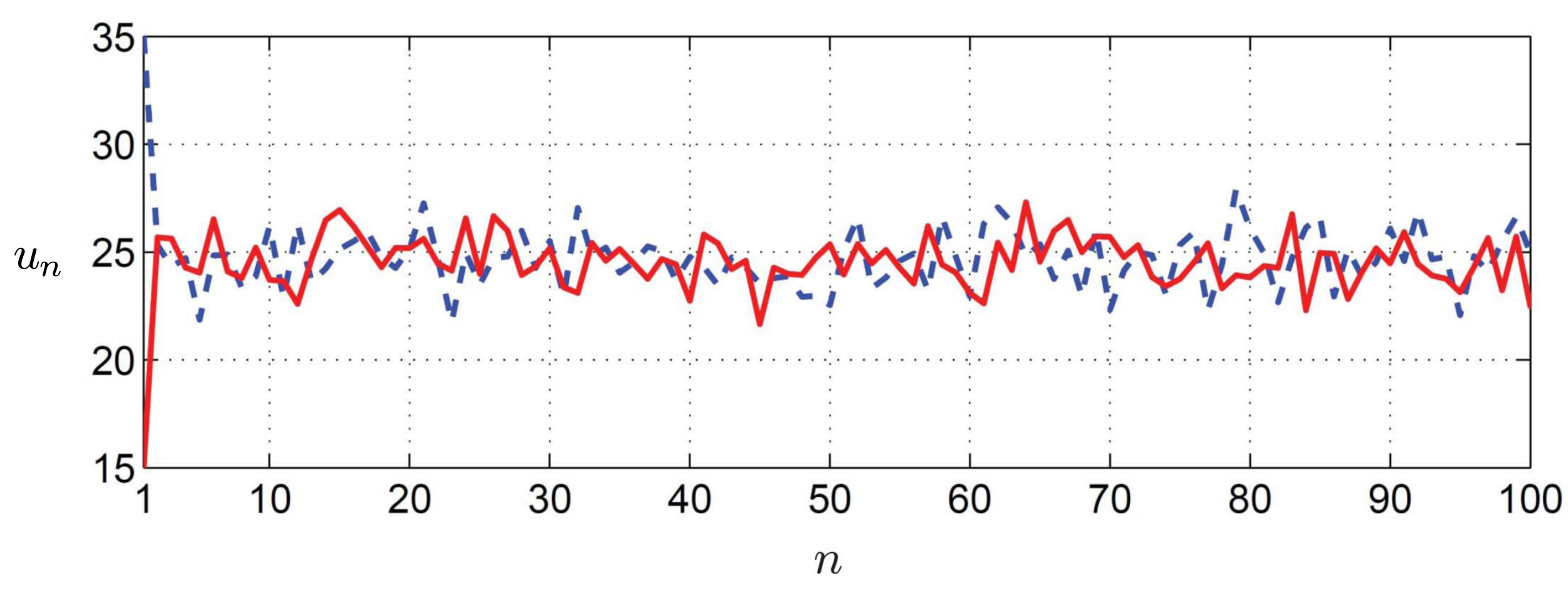}
{\small \caption{Petri net - threshold control.  $t_{f}=1,000$, $n=1,\ldots,100$}}
\end{figure}

\vspace{.1in}
\begin{figure}[h]
\centering
\includegraphics[width=0.65\textwidth]{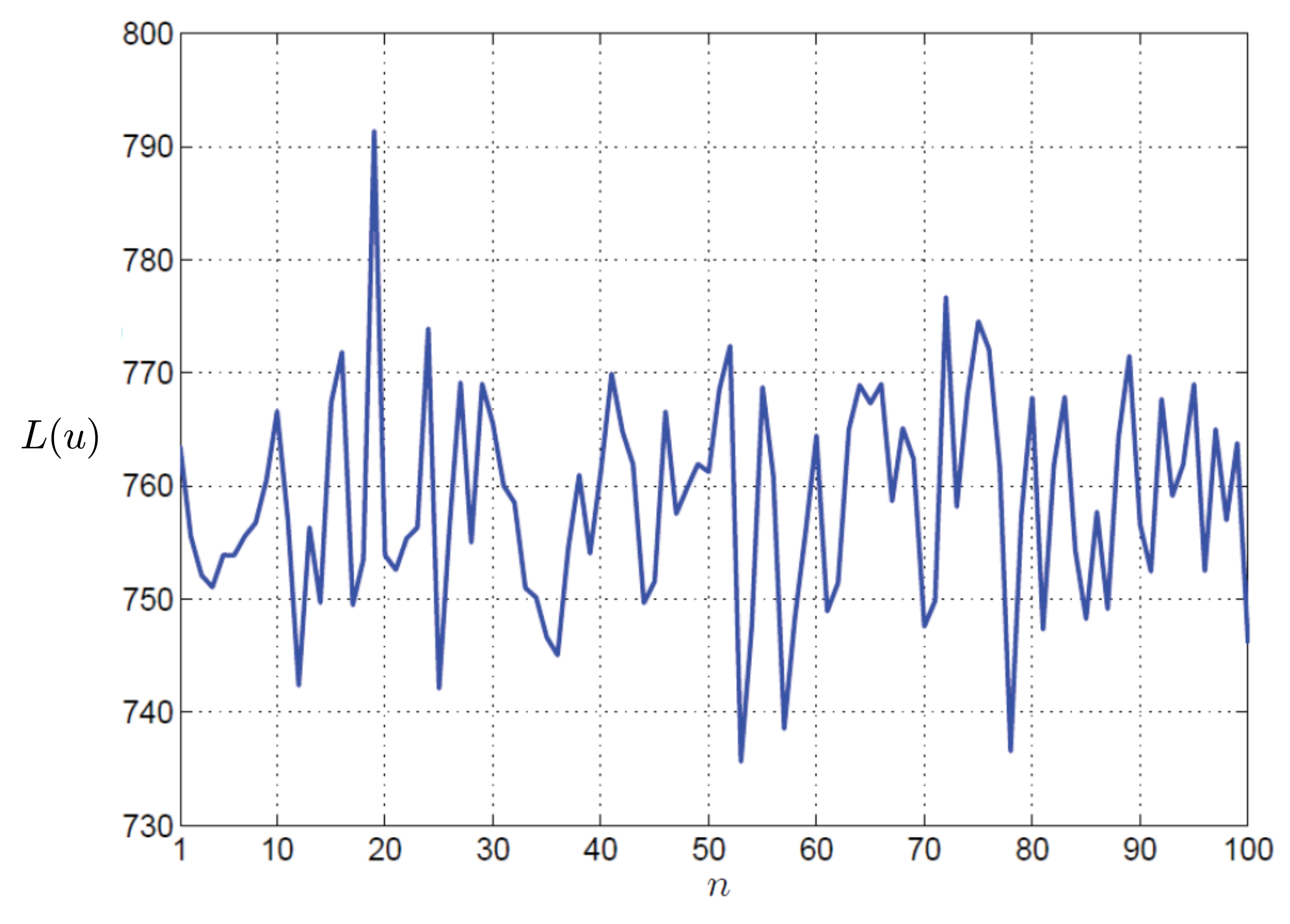}
{\small \caption{Variations in $L(u)$: 100 independent simulations}}
\end{figure}

\section{Application to throughput regulation in computer processors}\label{sec:chip}
The problem of maintaining stable instruction execution rates in computer
processors arises
in several application areas. For instance, in real-time
applications, guaranteeing a constant instruction execution rate
simplifies otherwise complex real-time scheduling tasks
by enabling tighter bounds on deadlines to be met \cite{Shakkottaia01}. Variable rate
instruction execution often leads to a reliance on worst case
execution-time estimates that far exceed the average execution
times, resulting in a  poorer processor utilization than is
necessary.  A second application concerns a tradeoff between
instruction execution rate and power dissipation.   Dynamic Voltage
Frequency Scaling (DVFS) can be used to navigate this tradeoff to
arrive at an instruction rate that is the most power-efficient~\cite{c3}\cite{c4}.  This approach, described in the
Introduction, was followed by~\cite{Almoosa12a} for throughput regulation in
processor cores. In this section we pursue a similar approach but use a
more-precise system-model for the throughput simulation and hence obtain
better results, including far-faster convergence, as described in the sequel.

We consider a multiprogrammed multi-core processor where programs are
assigned for execution to the cores by the operating system.  Each
core is assigned an instruction-execution rate setpoint by a
supervisory controller, and has to control (regulate) its instruction
rate to that level.  In this paper we are not concerned with the way
these setpoint levels are assigned, and consider them as given and
fixed.  Furthermore, we assume that each core is in a separate clock domain
and can independently control its own clock rate.  Each core exploits
instruction-level parallelism (ILP) utilizing an Out-of-Order (OOO)
technology whereby instructions may complete execution in an order different
from program order (i.e., out of order). OOO execution enables
instruction execution to be limited only by data dependencies rather than
the order in which they appear in the program thereby serving as the
primary means for exploiting ILP in modern processors.

A high-level functional and logical description of programs' execution in an OOO core is depicted in Figure 14.
Instructions are fetched sequentially from memory and placed in the Instruction queue. There they are processed by
functional units which can be thought of as servers in the queueing parlance.  It is often the case
that there are enough functional units available for concurrent execution of all of the instructions
in the queue. What holds up the execution of an instruction is its dependency on data that would become
available from the execution of another instruction. For example, the instruction of adding two variables, $x$ and $y$,
requires that both variables have numerical values. $x$ may be obtained from another numerical instruction,
while $y$ may be fetched from memory. Thus, it is evident that the addition of $x$ and $y$ cannot commence before the
instructions of computing $x$ and fetching $y$ are completed.

Instructions can be classified as {\it data dependent} vs. {\it data independent} according to whether
they depend on data provided by the execution of other instructions. If an instruction is data dependent then
its execution can commence as soon as all such data and a functional unit become available. If the instruction
is data independent then its execution can start as soon as a functional unit is available.

Memory instructions can take one-to-two orders of magnitude more time to execute than computational
instructions. Therefore most architectures make use of a hierarchical memory
arrangement where on-chip cache access takes less time than external memory such as DRAM. First the cache is searched,
and if the variable is found there then it is fetched and the instruction is completed. If the variable
is not stored in cache then a cache miss occurs and a cache line
(containing the variable) is fetched from external memory
(typically DRAM) and placed in the cache. The variable is then
accessed and the instruction is completed. We can
think of all external-memory (non-cache) accesses as placed in a
First-in-First-out queue designated by the {\it Memory} box in the
figure. Furthermore, this queue has a finite buffer, and when it becomes
full the entire memory access, including from the cache, is stalled.

During its processing the instruction is still stored in the Instruction queue, where it has been placed since its
arrival at the start of the aforementioned process. Although it may be executed concurrently with other instructions
subject to data-dependency constraints, it may not be released from the queue upon completion of its execution; in fact, its
release time is the later of its completion time and the release time of the previous instruction. However, the variables it computes
become available to other instructions upon its completion and not release. We point out that the term {\it release} used here
is called {\it commitment} in the parlance of computer architectures. It is this quantity that is used to compute the instruction's throughput, as
will be made clear in the sequel.

\begin{figure}[h]
\centering
\includegraphics[width=0.95\textwidth]{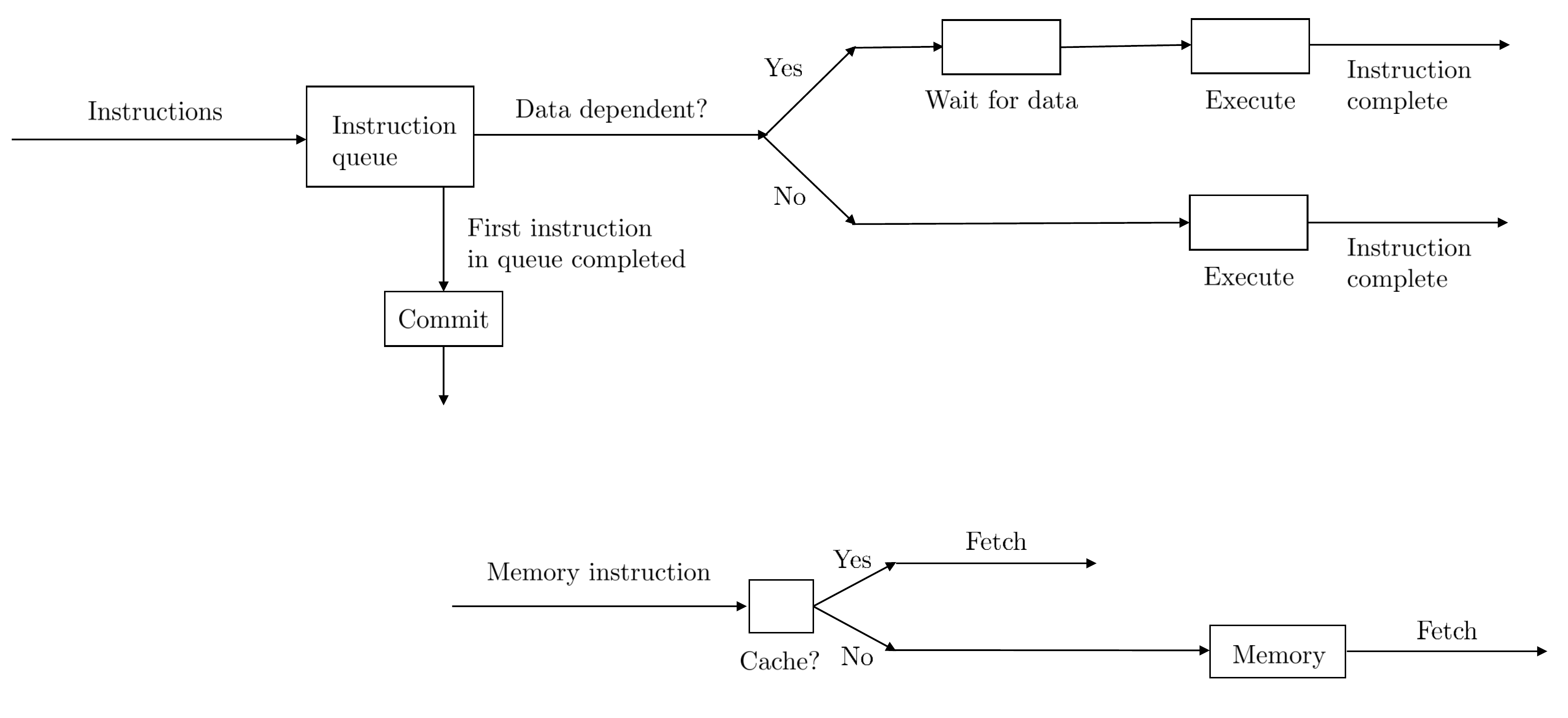}
{\small \caption{Functional description of OOO processing}}
\end{figure}

A high-level
description of the hardware  of an
OOO core is depicted in Figure 15 (see~\cite{c5} for a more detailed
description). The acronyms in the various blocks indicate the
following entities: IQ - Instruction Queue; ROB - Reorder Buffer; RS -
Reservation Station; FU - Functional Unit; RF - Register File;
L/S~Q - Load/Store Queue; MSHR - Miss Status
Handling Register; and MEM - memory other than cache, typically DRAM.
During a program's execution, the instruction unit of the processor fetches
program instructions sequentially from memory and places them in the
Instruction Queue (IQ). Instructions are issued from the IQ to a
reservation station waiting for access to the corresponding functional
unit. An instruction is granted access  (starts processing) after (i)  all its operands are available, for
example, the result from  preceding instructions, and (ii) the
functional unit is available. When an instruction
is issued to an RS, it is also allocated an entry in the ROB in the original program
order.

Instructions in the RS are issued to the functional units when their
operands are available (note that this issue order may be different from
program order hence the OOO designation). While the instruction is
waiting for its operands, it is said to be {\em stalled}. Instruction
execution results are broadcast to all reservation stations and the
ROB - they are stored along with the instruction's entry in the
ROB. When an instruction reaches the head of the ROB and it has
completed execution, its result is placed in the RF. This
process is called {\em instruction committment (IC)}. The {\it IC}
stage basically signifies the termination of the instruction processing. Note that
while instructions may complete execution out of order, they are
committed in program order. For the purpose of this discussion and its
related model, we can view the IQ and ROB as a single block where
instructions are buffered in order according to their issue (arrival)
order from the instruction fetch unit of the processor.

\begin{figure}[h]
\centering
\includegraphics[width=0.65\textwidth]{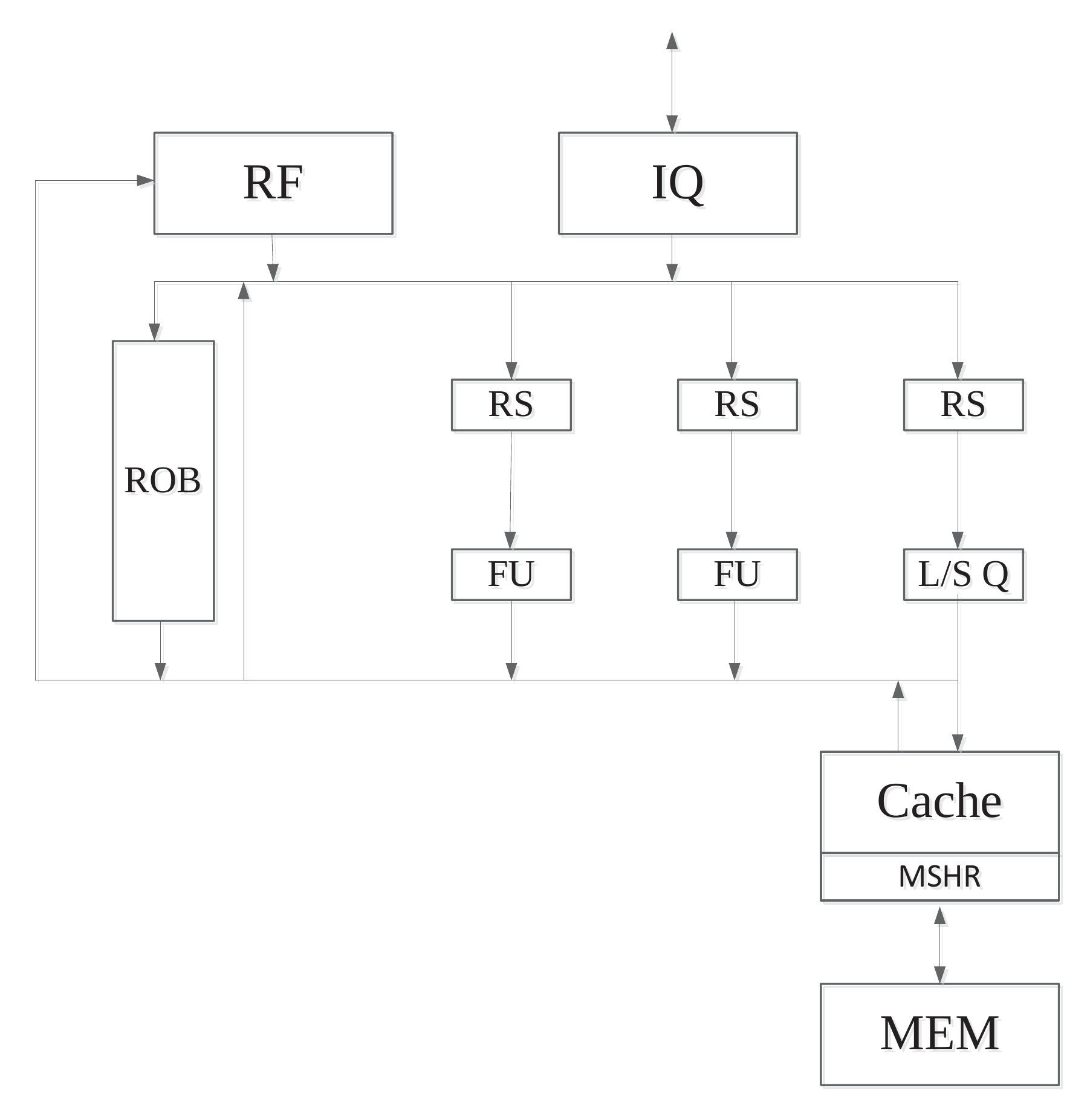}
{\small \caption{High-level OOO architecture}}
\end{figure}


{\it Memory} instructions such as load and store instructions are
directed to the Load/Store Queue (L/SQ) shown in Figure 15. Consider
the execution of a load instruction. First an attempt is made to
access the data from the data cache. If that is successful (cache
hit), the data is read from the cache line and sent to the
corresponding instruction entry in ROB.  If the requested data is not
in the cache (cache miss) it has to be fetched from the next level of
the memory hierarchy, for example main memory, typically DRAM.  In that case a request is sent to the Miss Status
Handling Register (MSHR), which serves as a finite-buffer queue for
buffering outstanding non-cache memory access requests.  When the MSHR
is full all new instructions that generate memory-accesses are
stalled.  Upon completion of the memory access, the data is sent to
the ROB and the instruction clears the MSHR.  Note that the MSHR is a
finite-buffer queue that holds only non-cache memory instructions. However,
when it becomes full, it halts all new memory access requests, and thus
 comprises a fairly nonstandard queueing model.

To quantify all of this (see also \cite{Chen15}), consider a control cycle comprised of
instructions $I_{i}$, $i=1,\ldots,M$ for a given $M>0$.  Let us set
the time to $t=0$ at the moment the first instruction is issued.  In
the framework of the closed-loop system defined by the regulation algorithm,
the control variable $u$ will be the clock frequency (rate) of
the core, but it is easier to carry out the IPA computations in terms
of the clock cycle, $\tau:=u^{-1}$.

For instruction $I_{i}$, define  $\xi_{i}$ to be its issue time in terms  of clock cycles (counting from the issue time of
$I_{1}$); $\alpha_{i}(\tau)$  - its execution starting time (in seconds), and
$\beta_{i}(\tau)$ - its completion time (in seconds), not to be confused with its commit
time which may occur later. Note that we use two kinds of timing variables, namely  in units of clock cycles or seconds. The  former can be measured  in real time for the purpose of the control law,
 while the latter
are used only in the analysis.

Consider first the case where   $I_{i}$ is not a memory instruction; memory instructions will be
handled later.
The issue time  of $I_{i}$ (in seconds) is $\xi_{i}\tau$, and we assume that there are available resources in the
ROB and RS so that the instruction is forwarded there at the same time,
$\xi_{i}\tau$.  If $I_{i}$ is data dependent,
denote by $k(i)$ the index (counter) of the instruction which is the last to provide a variable to $I_{i}$.
Then, we have that
\begin{equation}
\alpha_{i}(\tau)=\left\{
\begin{array}{ll}
\max\big\{\xi_{i}\tau,\beta_{k(i)}(\tau)\big\}+\tau, & {\rm if}\ I_{i}\ {\rm is\ data\ dependent}\\
\xi_{i}\tau+\tau, & {\rm otherwise}.
\end{array}
\right.
\end{equation}
As for the completion time, let
$\mu_{i}$ denote the execution time of $I_{i}$ in terms of clock cycles. Then,
\begin{equation}
\beta_{i}(\tau)=\alpha_{i}(\tau)+\mu_{i}\tau.
\end{equation}
Consider next the case where $I_{i}$ is a memory instruction. Upon its issuance at time $\xi_{i}\tau$ it is
directed to the L/SQ where a cache attempt is made. Let us regard the starting time of the cache attempt as
the starting time of the instruction's execution at the memory stage,
and  denote it by $\alpha_{i}(\tau)$. If the MSHR is full at time $\xi_{i}\tau$, let $\ell(i)$ denote the
index of the instruction at the head of the MSHR. Then,
\begin{equation}
\alpha_{i}(\tau)=\left\{
\begin{array}{ll}
\xi_{i}\tau+\tau, & {\rm if\ MSHR\ is\ not\ full\ at\ time}\ \xi_{i}\tau\\
\beta_{\ell(i)}(\tau)+\tau, & {\rm if\ MSHR\ is\ full\ at\ time}\ \xi_{i}\tau.
\end{array}
\right.
\end{equation}
Let $\nu_{i}$ denote the time (in units of clock cycles) it takes the L/SQ to process $I_{i}$
including the cache attempt. If the cache attempt result is successful (cache hit), then
\begin{equation}
\beta_{i}(\tau)=\alpha_{i}(\tau)+\nu_{i}\tau.
\end{equation}
On the other hand, in case of a cache miss, the instruction is directed to the MSHR for non-cache memory access.
Let $m_{i}$ denote the processing time of
 the instruction in the MSHR in units of clock cycles, and let
$MEM_{i}$ be the time it takes to access the memory.  Note that DRAM access is not governed by the core's clock  and hence it is not a function of $\tau$. Let $j(i)$ denote the index of the instruction prior to $I_{i}$   in the MSHR.
Then, the completion time
of $I_{i}$ is given by
\begin{equation}
\beta_{i}(\tau)=\max\big\{\alpha_{i}(\tau)+\nu_{i}\tau+m_{i}\tau+MEM_{i},\beta_{j(i)}(\tau)+\tau\big\}.
\end{equation}
Finally, the commit (departure) times of instructions $I_{i}$, $i=1,2,\ldots$, denoted by $d_{i}(\tau)$, are
given by the following recursive equation,
\begin{equation}
d_{i}(\tau)=\max\big\{\beta_{i}(\tau),d_{i-1}(\tau)\big\}+\tau.
\end{equation}

Now recall that the control cycle consists of $M$ instructions, and define the
average throughput by
\begin{equation}
y:=L(u)=\frac{M}{d_{M}(u)}.
\end{equation}
With $u=\tau^{-1}$, its IPA derivative  is
\begin{equation}
L^{\prime}(u)=\frac{1}{M}\Big(\frac{y}{u}\Big)^2d_{M}^{\prime}(\tau),
\end{equation}
where we assume  that the throughput $y=L(u)$ can be  computed from the system.
The IPA term $d_{M}^{\prime}(\tau)$ is computable  in a recursive manner  as follows. By Equations (27)--(32),
for every
$i=1,\ldots,M$: If $I_{i}$ is not a memory  instruction, then
\begin{equation}
\alpha_{i}^{\prime}(\tau)=\left\{\begin{array}{ll}
\beta_{k(i)}^{\prime}(\tau)+1, & {\rm if}\ I_{i}\ {\rm is\ stalled\ due\ to\ data\ dependency}\\
\xi_{i}+1, & {\rm otherwise},
\end{array}
\right.
\end{equation}
and
\begin{equation}
\beta_{i}^{\prime}(\tau)=\alpha_{i}^{\prime}(\tau)+\mu_{i}.
\end{equation}
On the other hand, if $I_{i}$ is a memory instruction, then
\begin{equation}
\alpha_{i}^{\prime}(\tau)=\left\{
\begin{array}{ll}
\beta_{\ell(i)}^{\prime}(\tau)+1, & {\rm if}\ I_{i}\ {\rm is\ stalled\ due\ to\ full\ MSHR}\\
\xi_{i}+1, & {\rm otherwise};
\end{array}
\right.
\end{equation}
as for $\beta_{i}^{\prime}(\tau)$, if $I_{i}$ results in a cache hit, then
\begin{equation}
\beta_{i}^{\prime}(\tau)=\alpha_{i}^{\prime}(\tau)+\nu_{i},
\end{equation}
and in the event of a cache miss,
\begin{equation}
\beta_{i}^{\prime}(\tau)=\left\{
\begin{array}{ll}
\alpha_{i}^{\prime}(\tau)
+\nu_{i}+m_{i}, & {\rm if}\ I_{i}\ {\rm is\ first\ in\ the\ MSHR\ queue\ by\ the\ time}\\
 & {\rm its\ variable\ is\ read\ from\ memory}\\
\beta_{j(i)}(\tau)+1, & {\rm otherwise}.
\end{array}
\right.
\end{equation}
It is reasonable to assume  that the quantities $\xi_{i}$, $\mu_{i}$, $\nu_{i}$, and $m_{i}$ can be read from the system in real time
during instruction execution,
 and hence the computation of $\beta_{i}^{\prime}(\tau)$, $i=1,\ldots,$ can be performed in
real time in
a recursive fashion via Equations (35)-(39).
Finally, by (32),
\begin{equation}
d_{i}^{\prime}(\tau)=\left\{
\begin{array}{ll}
d_{i-1}^{\prime}(\tau)+1, & {\rm if}\ I_{i}\ {\rm is\ stalled\ in\ ROB\ following\ its\ execution}\\
\beta_{i}^{\prime}(\tau)+1, & {\rm if}\ I_{i}\ {\rm is\ comitted\ right\ after\ its\ execution.}
\end{array}
\right.
\end{equation}
This yields $d_{M}^{\prime}(\tau)$, and hence  $L^{\prime}(u)$ via Equation (34).
We point out that this IPA derivative is biased due to the fact that
the DRAM and other non-cache memory accesses are not controlled by the
core's clock - the memory system is in a different clock
domain. We note that the  time required for such memory access  typically is one order
of magnitude longer than  a cache-access time and can be two orders of
magnitude longer than computing  instructions. Therefore
we expect the regulation technique to perform better when applied to
computation-intensive programs rather than to memory-intensive programs.  This is evident
from the simulation results which will be presented in the following paragraphs.

We implemented the control algorithm  using Manifold,
a cycle-level, full system discrete event simulation platform for
multi-core architectures~\cite{Yalamanchili}. A Manifold model boots a Linux
operating system and executes stock 32-bit x86 binaries.  Our
experiments were applied to two programs from the SPLASH-2 suite of
benchmark programs~\cite{c11}, {\it Barnes} and {\it
  Radiosity}. Barnes is compute intensive while Radiosity is memory
intensive. For both cases the Manifold processor model used is similar
to the Intel Nehalem micro-architecture comprised of four cores, each
having its own L1 cache and sharing an L2 cache \cite{Thomadakis11}. Each core is in a
separate clock domain that is independently controlled. At the start of
the simulation the application is emulated   for
a million instructions (out of an order of   $10^{12}$
instructions)
in
order to warm up the architecture state. At this point, cycle-level
timing simulation is begun over program regions of interest.

The control cycle for each core consists of $M=10,000$ instructions.
Thus, for a given core, the control variable during the $nth$ cycle is
$u_n$, and $y_n$ is the instruction throughput computed at the end of the cycle
via Eq. (33) (with the index $n$ added). The IPA  derivative is computed via Eq. (34).
It is not unbiased since the sample-performance function $L(u)$ is not necessarily continuous.
This is
largely due to the
 possibility of instruction stalls when the MSHR queue becomes full; see \cite{Ho91} for the relation between discontinuous
 sample performance functions and the biasedness of IPA. However, we believe that the error introduced by the bias generally is not
 large enough to prevent convergence of the regulation algorithm. Therefore, in contrast with the
 case of the finite-buffer queue discussed in Section 3.2, we do not resort to a fluid queueing model but rather compute the IPA derivatives directly from the discrete model according to Eqs. (34)~-~(40).

In the simulation experiments we took the target
 instruction rates (setpoints) for Cores 0-3 to be  1.0 Giga
Instruction Per Second (GIPS), 1.5 $GIPS$, 2.0 $GIPS$, and 2.5 $GIPS$,
respectively. Figure 16 shows simulation results for the benchmark
 {\it Barnes} executing at all four cores, and they indicate convergence-times of the
algorithm between 0.02 $ms$ (Core 3) and 0.08 $ms$ (Core
2).\footnote{Using a cruder model, Reference \cite{Almoosa12a} reports
  convergence in about 1 ms.}  The apparent oscillations after
convergence are due to variations in the programs' instruction loads, and their
magnitudes are within 10\% of the respective target values. However,
the average instruction rates from time 0.1 ms to the final time (0.25
ms) are 0.9998, 1.5025, 1.9997, and 2.4985, which are within 0.2 \% of the respective   target
values.

  \vspace{.2in}
\begin{figure}[h]
	\centering
	\includegraphics[width=1\textwidth]{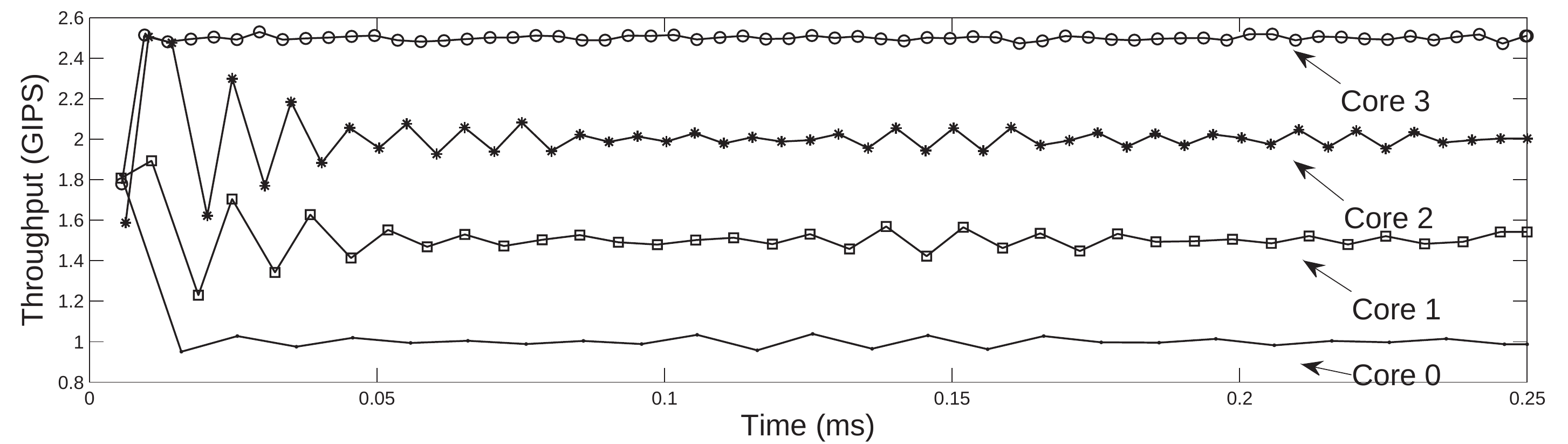}
	{\small \caption{Instruction rate regulation: {\it Barnes}}}
	\label{fig:fig12}
\end{figure}

For the Radiosity benchmark we set the target instruction rates  for Cores 0-3
to 1.0 $GIPS$, 1.3 $GIPS$, 1.5 $GIPS$, and 1.7 $GIPS$, respectively.
Typical results are shown in Figure 17, where we discern convergence
times between 0.12 $ms$ to 0.14 $ms$, with subsequent oscillations
below 15\% of the respective target values.  The slower convergence as
compared to Barnes is due to the fact that Radiosity is more memory
intensive.

\begin{figure}[h]
	\centering
	\includegraphics[width=1\textwidth]{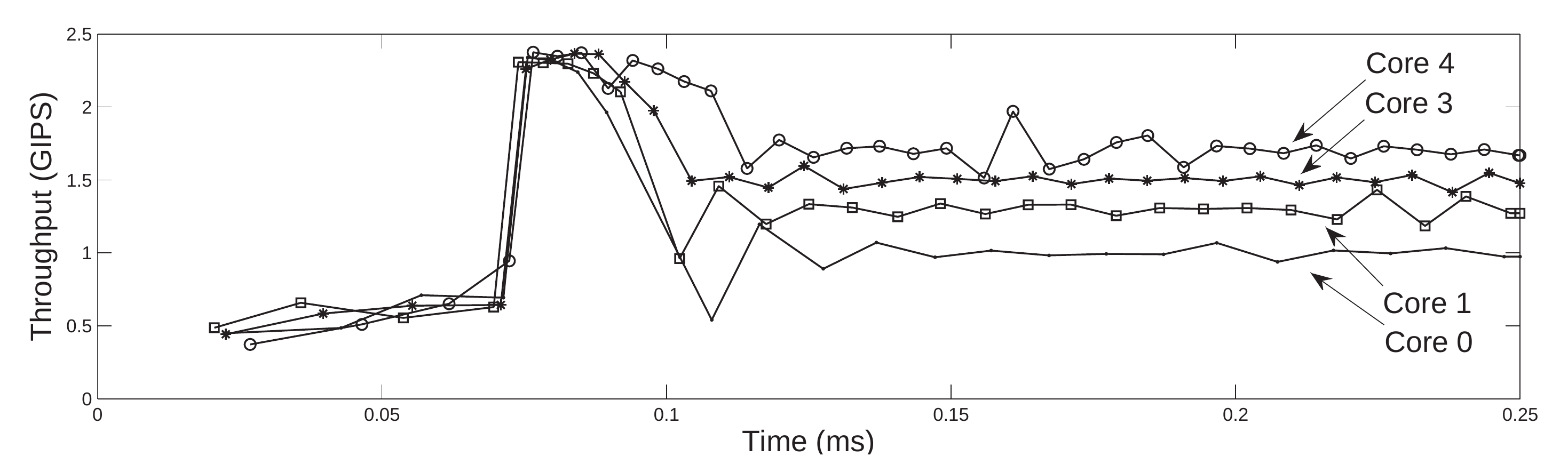}
	{\small \caption{Instruction rate regulation: {\it Radiosity}}}
	\label{fig:fig13}
\end{figure}

One way to reduce the oscillations is to scale the integrator's gain
in Equation (1) by a constant $k\in(0,1)$, thereby replacing (1) by
the following equation,
\begin{equation}
u_{n}=u_{n-1}+kA_{n}e_{n-1}.
\end{equation}
After some experimentation we chose $k=0.2$. This resulted  in  reductions
in the oscilations' magnitudes  from 10\% to 5\% for Barnes, and  from 15\% to 10\%
for radiosity. The results are shown in Figures 18 and
19, respectively.

\vspace{.2in}
\begin{figure}[h]
	\centering
	\includegraphics[width=1\textwidth]{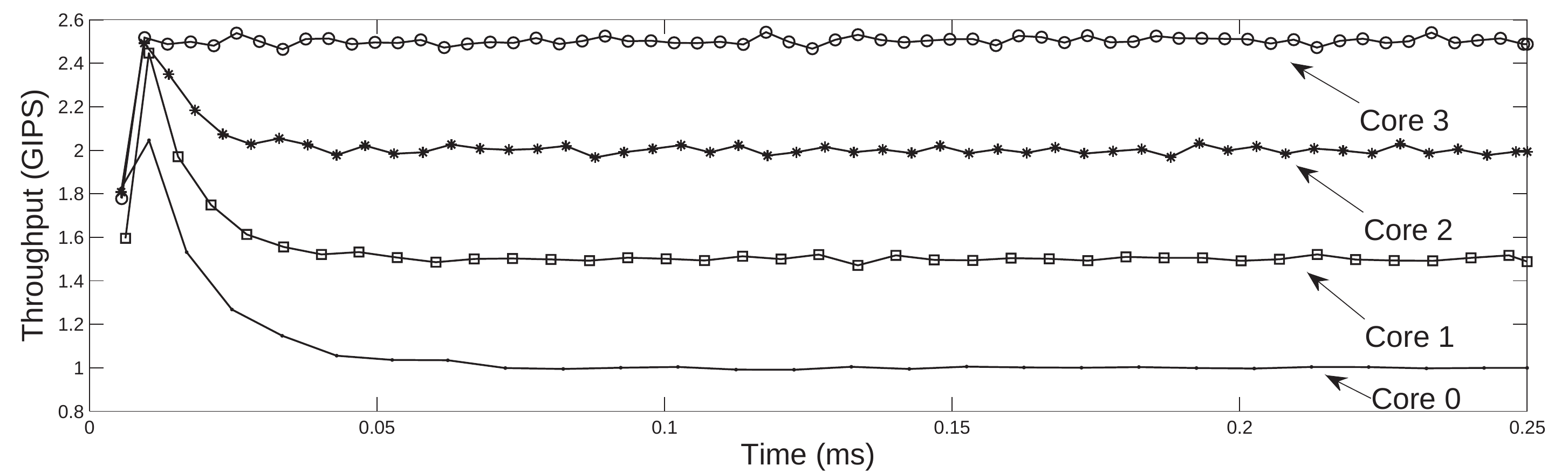}
	{\small \caption{Instruction rate regulation (modified algorithm): {\it Barnes}}}
	\label{fig:fig14}
\end{figure}

\vspace{.2in}
\begin{figure}[h]
	\centering
	\includegraphics[width=1\textwidth]{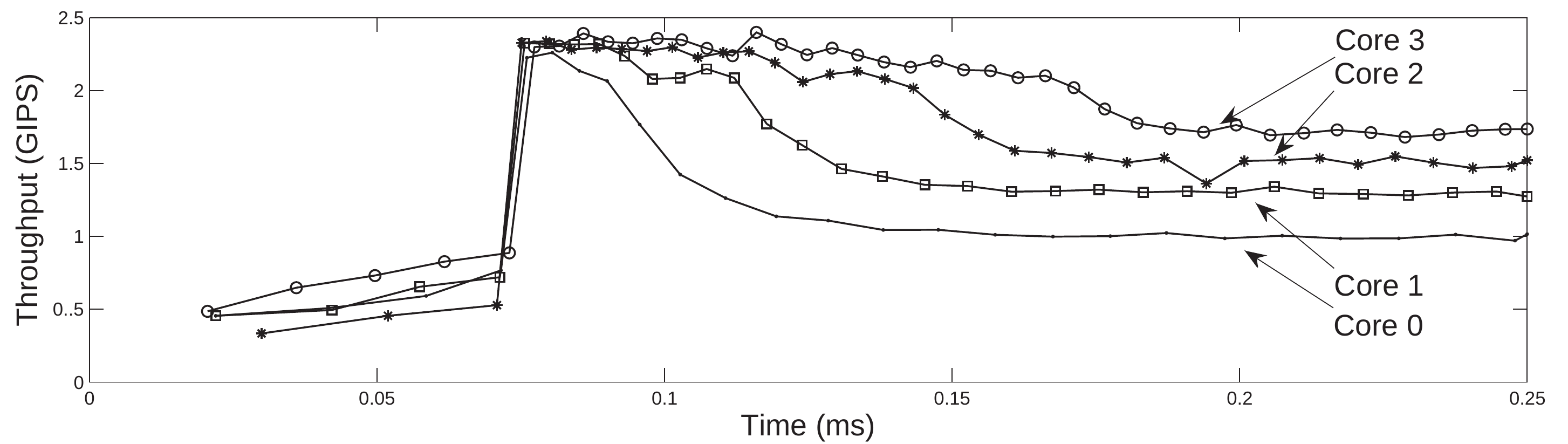}
	{\small \caption{Instruction rate regulation (modified algorithm): {\it Radiosity}}}
	\label{fig:fig15}
\end{figure}

\section{Conclusions}
This paper describes an IPA-based approach to performance regulation in stochastic timed discrete event dynamic systems.
The considered problem is to control the output of the system so as to have it track (asymptotically) a given
performance reference
in the face of variations in the system's characteristics. The proposed approach is based on
an integrator with adaptive gain. The system's plant is represented via a discrete-event or hybrid model, and the controller's gain is inverse proportional
to the IPA derivative of its plant function with respect to the control parameter.

The paper summarizes the regulation technique and presents several examples, which highlight the fact that it can
work  well despite  errors in estimating  the IPA gradient.  The examples include delay and loss in a single queue, inventory control in a Petri-net model of
a production system, and throughput regulation in a computer core. Extensions and future investigations will
focus on emerging problems in  various application  areas.

\section{Appendix}
This section provides proofs to Lemma 2.2 and Proposition 2.3.\\ \\
{\bf Proof of Lemma 2.2.} Consider a closed, finite-length interval $I$ where $g(\cdot)$ is monotone-nondecreasing and convex throughout $I$. We next prove the assertion of the lemma for this case, while the situations where
 $g(\cdot)$ is monotone-non-increasing or concave can be analyzed by similar arguments and hence their proofs are omitted.

 By the monotonicity of $g(\cdot)$ we have that
 $g^{\prime}(u)\geq 0~\forall u\in I$, and we will use this fact throughout the forthcoming analysis.
 By convexity of $g(\cdot)$  we have the following inequalities for every $u\in I$ and $\Delta u\in R$ such that
 $u+\Delta u\in I$:
 \begin{equation}
 g(u)+g^{\prime}(u)\Delta u\leq g(u+\Delta u)\leq g(u)+g^{\prime}(u+\Delta u)\Delta u.
 \end{equation}
 Given  $M\geq 1$, $\beta\in(0,1)$, and $\alpha\in(0,1)$, suppose that, for some $j=1,\ldots$,  $u_{j-1}\in I$,
 ${\cal E}_{j-1}<\alpha$, and  ${\cal G}_{j-1}<\beta$.
 We next derive upper bounds on the ratio $|g(u_{j}|/|g(u_{j-1})|$ (for the case where $g(u_{j-1})\neq 0$)  in terms of $M$, $\beta$, and $\alpha$. The
 analysis concerns four separate cases: (i) $g(u_{j-1})\leq 0$ and $g(u_{j})\leq 0$; (ii)  $g(u_{j-1})\leq 0$ and $g(u_{j})\geq 0$;
 (iii) $g(u_{j-1})\geq 0$ and $g(u_{j})\leq 0$; and (iv) $g(u_{j-1})\geq 0$ and $g(u_{j})\geq 0$.
 The results, to be derived in the following paragraph, are:
 \begin{itemize}
 \item
 Case (i):
 \begin{equation}
 |g(u_{j})|\leq\Big(1-\frac{1-\beta}{1+\alpha}\Big)|g(u_{j-1})|.
 \end{equation}
 \item
 Case (ii):
 \begin{equation}
 |g(u_{j})|\leq\Big(M\frac{1+\beta}{1-\alpha}-1\Big)|g(u_{j-1})|.
 \end{equation}
 \item
 Case (iii):
 \begin{equation}
 |g(u_{j})|\leq\Big(\frac{1+\beta}{1-\alpha}-1\Big)|g(u_{j-1})|.
 \end{equation}
 \item
 Case (iv):
 \begin{equation}
 |g(u_{j})|\leq\Big(1-\frac{1}{M}\frac{1-\beta}{1+\alpha}\Big)|g(u_{j-1})|.
 \end{equation}
 \end{itemize}
 We next prove these inequalities. In all cases we consider Eq. (42) with $u=u_{j-1}$ and
 \[
 \Delta u=-\frac{1}{g^{\prime}(u_{j-1})+\phi_{j-1}}\big(g(u_{j-1})+\psi_{j-1}\big).
 \]
 By Eq. (10), $u+\Delta u=u_j$.

 {\it Case (i): $g(u_{j-1})\leq 0$ and $g(u_{j})\leq 0$}. The left inequality of Eq. (42)  means that
 \begin{equation}
 g(u_{j-1})-g^{\prime}(u_{j-1})\frac{1}{g^{\prime}(u_{j-1})+\phi_{j-1}}\big(g(u_{j-1})+\psi_{j-1}\big)\leq g(u_j),
 \end{equation}
 and with a straightforward algebra,
 \begin{equation}
 g(u_{j})\geq g(u_{j-1})-\frac{1}{1+\frac{\phi_{j-1}}{g^{\prime}(u_{j-1})}}\Big(1+\frac{\psi_{j-1}}{g(u_{j-1})}\Big)g(u_{j-1}).
 \end{equation}
 By definition of ${\cal E}_{j-1}$ and ${\cal G}_{j-1}$, $|\frac{\phi_{j-1}}{g^{\prime}(u_{j-1})}|\leq{\cal E}_{j-1}$ and
 $|\frac{\psi_{j-1}}{g(u_{j-1})}|\leq{\cal G}_{j-1}$; and by assumption, ${\cal E}_{j-1}\leq\alpha$ and
 ${\cal G}_{j-1}\leq\beta$; hence, and by (48) and the assumption that $g(u_{j-1})\leq 0$, we have that
 \begin{equation}
 g(u_{j})\geq\Big(1-\frac{1-\beta}{1+\alpha}\Big)g(u_{j-1}).
 \end{equation}
 By the assumption that $g(u_j)\leq 0$, Eq. (43) follows.

 {\it Case (ii): $g(u_{j-1})\leq 0$ and $g(u_{j})\geq 0$}. The right inequality of (42) means that
 \begin{equation}
 g(u_{j})\leq g(u_{j-1})-\frac{g^{\prime}(u_j)}{g^{\prime}(u_{j-1})+\phi_{j-1}}(g(u_{j-1})+\psi_{j-1}),
 \end{equation}
 and multiplying and dividing the RHS of (50) by $g^{\prime}(u_{j-1})$ and $g(u_{j-1})$ we obtain that
 \begin{equation}
 g(u_{j})\leq g(u_{j-1})-\frac{g^{\prime}(u_j)}{g^{\prime}(u_{j-1})}\frac{1}{1+\frac{\phi_{j-1}}{g^{\prime}(u_{j-1})}}
 \Big(1+\frac{\psi_{j-1}}{g(u_{j-1})}\Big)g(u_{j-1}).
 \end{equation}
By definition of ${\cal E}_{j-1}$ and ${\cal G}_{j-1}$,
$|\frac{\phi_{j-1}}{g^{\prime}(u_{j-1})}|\leq{\cal E}_{j-1}$ and
 $|\frac{\psi_{j-1}}{g(u_{j-1})}|\leq{\cal G}_{j-1}$; and by assumption, ${\cal E}_{j-1}\leq\alpha$,
 ${\cal G}_{j-1}\leq\beta$, and
 $\frac{g^{\prime}(u_{j})}{g^{\prime}(u_{j-1})}\leq M$; hence, and by (51) and the assumption that $g(u_{j-1})\leq 0$, we have that
 \begin{equation}
 g(u_{j})\leq\Big(1-M\frac{1+\beta}{1-\alpha}\Big)g(u_{j-1}).
\end{equation}
 Since by assumption $g(u_{j-1})\leq 0$ and $g(u_{j})\geq 0$, Eq. (44) follows.

 {\it Case (iii): $g(u_{j-1})\geq 0$ and $g(u_{j})\leq 0$}. As in Case (i), Eq. (48) is satisfied, and since $g(u_{j-1})\geq 0$,
 \begin{equation}
 g(u_{j})\geq g(u_{j-1})-\frac{1+\beta}{1-\alpha}g(u_{j-1})=\Big(1-\frac{1+\beta}{1-\alpha}\Big)g(u_{j-1}).
 \end{equation}
 Since $g(u_{j})\leq 0$, Eq. (45) follows.

 {\it Case (iv):  $g(u_{j-1})\geq 0$ and $g(u_{j})\geq 0$}.
 As in the analysis for Case (ii), Eq. (51) applies, and by the definition of $M$,
 $\frac{g^{\prime}(u_{j})}{g^{\prime}(u_{j-1})}\geq\frac{1}{M}$. Therefore
 \begin{equation}
 g(u_{j})\leq\Big(1-\frac{1}{M}\frac{1-\beta}{1+\alpha}\Big)g(u_{j-1}),
 \end{equation}
 which is Eq. (46).

 Fix $M>1$ and $\beta\in(0,M^{-1})$. Consider $\alpha\in(0,1)$, and $n=1,\ldots$, and suppose that all of the conditions specified
 in the assertion of the lemma are satisfied. We next prove that Eq. (11) is satisfied for an $\alpha>0$ and a corresponding
 $\theta\in(0,1)$. There are two scenarios to consider: (a) $m_{n}=n+1$, and (b) $m_{n}>n+1$.

 Scenario (a) means that $g(u_{n})g(u_{n-1})\geq 0$ and arises when either Case (i) or Case (iv) occur. In Case (i), Eq.
 (43) is in force, and since $\big(1-\frac{1-\beta}{1+\alpha}\big)\big|_{\alpha=0}=\beta<1$, it follows that there exists
 $\alpha_{1}\in(0,1)$ and $\theta_{1}\in(0,1)$ such that, for every $\alpha\in(0,\alpha_{1})$, Eq. (11) is satisfied with $\theta_{1}$ in
 lieu of $\theta$. Similarly in Case (iv) and Eq. (46); since $\big(1-\frac{1}{M}\frac{1-\beta}{1+\alpha}\big)\big|_{\alpha=0}=1-\frac{\beta}{M}<1$,
 there exists $\alpha_{4}\in(0,1)$ and $\theta_{4}\in(0,1)$ such that, for every $\alpha\in(0,\alpha_{4})$, Eq. (11)
 is satisfied with $\theta_{4}$ in lieu of $\theta$.

 Scenario (b) corresponds to either Case (ii) or Case (iii). In Case (iii), where $g(u_{n-1})\geq 0$ while $g(u_{n})\leq 0$,
 Eq. (45) is satisfied. Note that $\big(\frac{1+\beta}{1-\alpha}-1\big)\big|_{\alpha=0}=\beta<1$, and therefore, there exists
 $\alpha_{3}\in(0,1)$ and $\theta_{3}\in(0,1)$ such that, if $\alpha\in(0,\alpha_{3})$, then
 \begin{equation}
 |g(u_{n})|\leq\theta_{3}|g(u_{n-1})|.
 \end{equation}
 In Case (ii), where $g(u_{n-1})\leq 0$ while $g(u_{n})\geq 0$, Eq. (44) holds, but unlike the other three cases,
 it is not true that $\big(M\frac{1+\beta}{1-\alpha}-1\big)\big|_{\alpha=0}<1$. A different argument is needed.

 Suppose first that Case (ii) holds at $u_{n-1}$,  Case (iii) will be considered later. By definition of $m_{n}$,  Case (iii) holds at $m_{n}-1$
  while Case (iv) is satisfied for all $j=n,\ldots,m_{n}-2$.  By Equations (44)-(46),
 \begin{equation}
 |g(u_{n_{n-1}})|\leq\Big(M\frac{1+\beta}{1-\alpha}-1\Big)\Big(\frac{1+\beta}{1-\alpha}-
 1\Big)\Big(1-\frac{1}{M}\frac{1-\beta}{1+\alpha}\Big)^{m_{n-1}-(n-1}|g(u_{n-1})|.
 \end{equation}
 The above analysis of Case (iv) showed that if $\alpha<\alpha_4$ then $\Big(1-\frac{1}{M}\frac{1-\beta}{1+\alpha}\Big)<1$, and therefore,
 \begin{equation}
 |g(u_{n_{n-1}})|\leq\Big(M\frac{1+\beta}{1-\alpha}-1\Big)\Big(\frac{1+\beta}{1-\alpha}-
 1\Big)|g(u_{n-1})|.
 \end{equation}
 But $\Big(M\frac{1+\beta}{1-\alpha}-1\Big)\Big(\frac{1+\beta}{1-\alpha}-
 1\Big)\Big|_{\alpha=0}=\big(M(1+\beta)-1\big)\beta<M\beta<1$, where the last two inequalities are due to the
  assumption that $M\beta<1$. Therefore, there exist $\alpha_{2}\in(0,1)$ and $\theta_{2}\in(0,1)$ such that, if
  $\alpha<\alpha_{2}$, Eq. (11) is satisfied with $\theta_{2}$ in lieu of $\theta$.

  Finally, in Case (iii) at $u_{n-1}$, Eq. (57) is provable in the same way as for Case (ii), and the conclusion is derivable in
  the same way as well.

  Now by defining $\alpha=\min\{\alpha_{i}:i=1,\dots,4\}$, Eq. (11) is satisfied with
  $\theta:=\max\{\theta_{i}:i=1,\ldots,4\}$. This completes the proof. \hfill $\Box$\\ \\
  {\bf Proof of Proposition 2.3.}
  Given $\eta>0$, $M>1$, and $\varepsilon>0$. Fix $\beta\in(0,M^{-1})$. Choose $\alpha\in(0,1)$ and $\theta\in(0,1)$ according
  to Lemma 2.2. In particular, as in the proof of Lemma 2.2 we can assume, by reducing $\alpha$ is necessary, that
  \begin{equation}
  \Big(M\frac{1+\beta}{1-\alpha}-1\Big)\Big(\frac{1+\beta}{1-\alpha}-1\Big)<1.
  \end{equation}
  Consider a closed, finite-length interval $I$ and a sequence of
  points $\{u_{n}\}_{n=1}^{\infty}$ satisfying the conditions of the proposition. Suppose, without loss of generality,
  that $g(\cdot)$ is monotone nondecreasing and convex
  on $I$. Since $|g^{\prime}(u)|\geq\eta$ for
  every $u\in I$, it follows by Eq. (10) that there exists $K>0$ (independent of $I$ or the sequence
  $\{u_{n}\}$) such that, for every $n=1,\ldots$,
  \begin{equation}
  |g(u_{n})|\leq K|g(u_{n-1})|.
  \end{equation}
  Fix $\varepsilon^{\prime}>0$ such that
  \begin{equation}
  2K\varepsilon^{\prime}<\varepsilon.
  \end{equation}
  Fix $\delta>0$ such that
  \begin{equation}
  \delta<\min\big\{\beta\varepsilon^{\prime},\frac{\varepsilon}{2}\big\}.
  \end{equation}
  By assumption ${\cal E}_{n}<\alpha$ and $|\psi_{n}|<\delta$ for all $n=1,\ldots,$. As a result of the inequality
  $|\psi_{n}|<\delta$,
  and by Eq. (61), if ${\cal G}_{n}:=\frac{|\psi_{n}|}{|g(u_{n})|}>\beta$ then $|g(u_{n})|<\varepsilon^{\prime}$.
  Thus, if $|g(u_{n})|>\varepsilon^{\prime}$ then Eqs. (43)-(46) hold with $n+1$ in lieu of $j$, and if
  $|g(u_{n})|<\varepsilon^{\prime}$ then
  \begin{equation}
  |g(u_{n+1})|\leq K\varepsilon^{\prime}.
  \end{equation}

Now consider a point $u_{n}$ such that $|g(u_{n})|>\varepsilon^{\prime}$ and hence
${\cal G}_{n}\leq\beta$.
Consider the four cases (i) - (iv) in the proof of Lemma 2.2.
In Case (i) and Case (iv), $m_{n}=n+1$, and Eq. (11) implies that
$|g(u_{n+1})|<|g(u_{n})|$. In Case (iii), Eq. (55) implies the same inequality, namely $|g(u_{n+1})|<|g(u_{n})|$.
Only in Case (ii) the reverse inequality is possible, namely that
$|g(u_{n+1})|\geq|g(u_{n})|$.

Suppose that Case (ii) holds at $u_{n}$. Then (by definition of Case (ii)) $g(u_{n})\leq 0$, $g(u_{j})\geq 0$ $\forall~j=n+1,\ldots,m_{n}-1$,
and $g(u_{m_{n}})\leq 0$. Moreover, for every $j=n+1,\ldots,m_{n}-1$, either Case (iii) or (iv) holds, and
therefore, either $|g(u_{j+1})|<|g(u_{j})|$ if $|g(u_{j})|\geq \varepsilon^{\prime}$, or
$|g(u_{j+1})|\leq K\varepsilon^{\prime}$ if
$|(g(u_{j})|\leq\varepsilon^{\prime}$. Consequently, and by Eq. (62), we have that
\begin{equation}
|g(u_{m_{n}})|\leq\max\big\{K\varepsilon^{\prime},|g(u_{n+1})|\big\}.
\end{equation}
As a result we have the following situation: If $|g(u_{j})|\leq\varepsilon^{\prime}$
for some $j=n,\ldots,m_{n}-1$, then
$|g(u_{m_{n}})|\leq K\varepsilon^{\prime}$. On the other hand,
if $|g(u_{j})|\geq\varepsilon^{\prime}$ for every $j=n,\ldots,m_{n}-1$,
then (by Lemma 2.2), Eq. (11) is in force.
This, in conjunction with Eq. (61), implies that
\begin{equation}
\limsup_{n\rightarrow\infty}|g(u_{n})|\leq K\varepsilon^{\prime}\leq\frac{\varepsilon}{2},
\end{equation}
 hence the inequality in Eq. (13).

Finally, Eq. (14) follows from (13) and the assumption that $\delta<\frac{\varepsilon}{2}$, as specified in Eq. (61).
This completes the proof. \hfill $\Box$

\end{document}